\documentclass[mathpazo]{cicp}
\usepackage{mathrsfs}
\usepackage{subcaption}
\usepackage{amssymb}
\usepackage{float}
\usepackage{amsbsy}
\usepackage[mathscr]{euscript}
\usepackage{bm}
\usepackage{csquotes}
\usepackage{microtype}
\usepackage{amsmath,amsthm,amssymb}
\usepackage{stmaryrd}
\usepackage{mathrsfs}
\usepackage{amsmath}
\usepackage{mathrsfs}
\usepackage{subcaption}
\usepackage{amssymb}
\usepackage{amsbsy}
\usepackage{hyperref}
\usepackage{cleveref}
\usepackage{graphicx}
\usepackage[mathscr]{euscript}
\usepackage{comment}
\usepackage[T1]{fontenc}
\usepackage{color}

\numberwithin{equation}{section} 

\begin{document}

%\journal{Computer Physics Communications}

\title{Discontinuous  collocation and symmetric integration methods for   distributionally-sourced hyperboloidal partial differential equations }

%%%%% author(s) :
\author[O'Boyle, M. et.~al.]{
Michael F. O'Boyle\affil{1}\comma\corrauth and
Charalampos Markakis\affil{2,3}}
\address{
\affilnum{1}\ Grainger College of Engineering, University of Illinois at Urbana-Champaign, Urbana, Illinois 61801, USA\\
\affilnum{2}\ Mathematical Sciences, Queen Mary University of London, E1 4NS, London, UK\\
\affilnum{3}\ NCSA, University of Illinois at Urbana-Champaign, Urbana, Illinois 61801, USA\\
}
\email{{\tt moboyle2@illinois.edu}}

\begin{abstract}
This work outlines a time-domain numerical integration technique for linear hyperbolic partial differential equations sourced by distributions (Dirac $\delta$-functions and their derivatives). Such problems arise when studying binary black hole systems in the extreme mass ratio limit.
We demonstrate that such source terms may be converted to effective domain-wide sources when discretized, and we introduce a class of time-steppers that directly account for these discontinuities in time integration.
Moreover, our time-steppers are constructed to respect time reversal symmetry, a property that has been connected to conservation of physical quantities like energy and momentum in numerical simulations.
To illustrate the utility of our method, we numerically study a distributionally-sourced wave equation that shares many features with the equations governing linear perturbations to black holes sourced by a point mass.
\end{abstract}

\ams{35Q75, 65D30, 65M20, 65M22, 65M70, 83-10, 83C25}
\keywords{Symmetric integration, Hermite integration, discontinuous collocation methods, black hole perturbation theory, hyperboloidal slicing}

\maketitle

\section{Introduction}
\label{sec:intro}
For radiation-reaction forces in electrodynamics and gravitational self-force corrections in general relativity \cite{Barack2009,Poisson2011}, the source of reduced field equations is distributional, leading to non-analyticity in the field at the particle's location. Previous solutions have ranged from multi-domain pseudospectral methods and discontinuous Galerkin methods with time-dependent mapping \cite{Canizares2009,Canizares2011a,Canizares2010c,Jaramillo2011,Field2009,Field:2010xn,Heffernan:2017cad}, to finite-difference methods based on null coordinates or representations of the Dirac $\delta$ function \cite{Barack:2000zq,Lousto:1997wf,Lousto:1999za,Barack:2002ku,Barack:2005nr,Haas:2007kz,Haas:2011np,Martel:2001yf,Sundararajan2007,Sundararajan2008,Harms2013}. The former two methods employ a ``Lagrangian'' perspective in that the particle is treated with co-moving coordinates (albeit in a decomposed domain), while the latter three employ an ``Eulerian'' perspective in that the particle moves with respect to fixed coordinates. The methods proposed here leverage the advantages of both perspectives, integrating high order derivative jumps for finite-difference or pseudospectral implementations in Eulerian coordinates without domain decomposition.

Existing methods for handling non-analyticities often incur significant computational cost and complexity, or sacrifice accuracy \cite{Boyd2001}. Post-processing of oscillatory data can recover spectrally convergent non-oscillatory solutions, but these techniques can also be computationally expensive and prone to complications \cite{Gottlieb1992,Gottlieb1994,Gottlieb1996,Gottlieb1995,Gottlieb1995a,Gottlieb1997,GottliebDavidSigal2003,2019JCoPh.390..527P,Jung2004,Archibald2003,Archibald2002,Archibald2004,Shu1995,Boyd2005}.

Existing methods for approximating piecewise smooth functions, such as those introduced by Krylov \cite{Krylov1907}, Lanczos \cite{Lanczos1966} and Eckhoff \cite{eckhoff1994,eckhoff1995,eckhoff1996,eckhoff1997,eckhoff1998,eckhoff2002}, or those that use polynomial correction terms as in Lipman and Levin \cite{Lipman2009}, often scale poorly or face other complications. In contrast, Lipman and Levin's method, uses moving least squares to determine the  location and magnitude of a discontinuity, and its computational cost scales more favorably, comparable to $O(N)$ rather than $O(N^2)$.

Explicit time-evolution schemes, such as classical Runge-Kutta methods, are prevalent in numerical relativity, despite their drawbacks: conditional stability (CFL limitation) and energy/symplectic structure violations in Hamiltonian systems. This ambiguity in energy loss cause can pose issues in gravitational wave (GW) computations. An alternative is a \textit{geometric integrator}, which upholds key aspects of Hamiltonian dynamics, like symplecticity or time-reversal symmetry \cite{hairer_geometric_2006,sanz-serna_numerical_1994}. While no method can universally preserve both energy and symplectic structure in a Hamiltonian system, this can be possible for quadratic Hamiltonians \cite{butcherSymmetricGeneralLinear2016,chartierAlgebraicApproachInvariant2006,zhongLiePoissonHamiltonJacobiTheory1988}, such as those arising in black hole perturbation theory.

Integrators which respect time-reversal symmetry are especially appealing to physics problems because they have been shown to have intimate connections to conservation in numerical dynamics \cite{moboylethesis,Markakis:2014nja,markakis2019timesymmetry}. In addition, these schemes are implicit, which endows them with unconditional stability \cite{brownMultiderivativeNumericalMethods1975,brownCharacteristicsImplicitMultistep1977}. By removing the CFL limit, these methods allow larger timesteps, overcome numerical instability, and speed up evolutions, making them suitable for long-time simulations of events like LISA band GWs.

In what follows, we consider distributionally-sourced partial differential equations (PDEs) of the 1+1 form  
\begin{equation}\label{eq:DistrForcedPDE}
    \Box \Psi_{\ell m} =  F_{\ell m}(t,x) ~ \delta'( x - \xi( t )) +G_{\ell m}(t,x) ~ \delta( x - \xi( t ))
\end{equation}
where $\Box $ is a generalized d`Alembert operator in a curvilinear coordinate chart $\{t,x \}$ and $x=\xi(t)$ is the worldline of a point particle in this chart. Here, we assume that the field $\Psi_{\ell m}$ has been expanded in spherical harmonic modes $\ell,m$ so that the remaining d'Alembert operator is 1+1. The main difficulty is that 
$\Psi_{\ell m}$ and its time or space derivatives are discontinuous accross the particle's wordline. The methods discussed here exploit the fact that the location and magnitude of these jump discontinuities can be determined in advance, in terms of 
$ F_{\ell m},  G_{\ell m}$ and $\xi(t)$. Below, for completeness, we
summarize earlier work on discontinuous collocation methods \cite{Markakis:2014nja} for solving equations of the above form. We then extend our discontinuous time-symetric integration formulae  from second to fourth order. We demonstrate the utility of these schemes in the case of a wave equation sourced by a moving particle on  hyperboloidal slices of Minkowski spacetime. Throughout this work, we adopt so-called geometric units where $c=1$ and time coordinates have the same units as space coordinates.

\section{Discontinuous collocation methods }\label{sec:DiscColl}

Here, for completeness, we outline the derivation of discontinuous collocation methods (\textit{cf.} \cite{Lipman2009,Markakis:2014nja,moboylethesis} for original details) which are used in \S4 to solve distributionally-sourced PDEs. The methods of \S2 are designed to model discontinuities in space, which suffices for static sources. For moving sources, one must also model discontinuities in time, as discussed in \S3.

\subsection{Piecewise smooth interpolation} \label{InterpolationSection}
Let $f:[a,b]\rightarrow \mathbb{R}$  be a $C^{k}$ function,  with the values $f_i=f(x_i)$ known at the ordered, distinct nodes  $x_i$ ($a\le x_0<x_1<\dots<x_N \le b$). Interpolation with a Lagrange polynomial of degree $N$ is limited by the degree of differentiability, that is, $N+1 \le k$. A degree of differentiability $k$ lower than $N+1$ necessitates  modification   of standard Lagrange interpolation. If  $f$ is a piecewise-$C^{N+1}$ function and  its jump discontinuities are known, then a simple generalization of the Lagrange method to any $k \ge -1$ can be found.

Let  the discontinuity be located at some $\xi \in (a,b)$ and  the  jumps in $f$ and its derivatives be given by:
\begin{equation} \label{fmJconditions2}
    {f^{(m)}}(\xi^+  ) - {f^{(m)}}(\xi^- ) := {J_m} < \infty ,\quad m = 0,1, \ldots ,\infty,
\end{equation}
where $\xi^-$ and $\xi^+$ denote the  limits to $\xi$ from below and above respectively. Approximate $f$  by a piecewise polynomial that interpolates the given $N+1$ nodes and  has the above specified derivative jumps at the discontinuity:
\begin{equation} \label{PolyN2DiscInterpol}
    p(x) =\theta (x - \xi ){p_ + }(x) + \theta (\xi  - x){p_ - }(x)
\end{equation}
where
\begin{equation} \label{PolyN1DiscInterpol}
    \theta (x)=\left\{ {\begin{array}{*{20}{c}} 1,&{x > 0}\\
    {1/2,}&{x = 0}\\
    0,&{x < 0} \end{array}} \right.
\end{equation}
is the Heaviside step function, and 
\begin{equation} \label{PolyN2}
    {p_ - }(x) = \sum\limits_{j = 0}^{N} {c_j^-(\xi) {x^j}},\quad{p_ + }(x) = \sum\limits_{j = 0}^{N} {c_j^+(\xi) {x^j}} 
\end{equation}
denote the left and right interpolating polynomials.
The polynomials $p_-,p_+$  and the piecewise polynomial  $p$ depend   on the location and  magnitude of the discontinuity.
To determine $c_j^\pm$ we use the method of undetermined coefficients. 
Given the ansatz \eqref{PolyN2DiscInterpol}, the collocation conditions become:
\begin{equation} \label{PolyNConds2}
 p(x_i)=f_i\Longleftrightarrow\   \left\{ \begin{array}{*{20}{c}}
    p_{+} (x_i) = f_i,&x_i > \xi\\
    p_{-} (x_i) = f_i,&x_i < \xi
    \end{array} \right..
\end{equation}The above $N+1$ collocation conditions only determine half of the ($2N+2$)   polynomial coefficients $c_j^\pm$. To close the system and determine  the other half, one  can impose  the first $N+1$ of the jump conditions  \eqref{fmJconditions2}. However, as higher order jumps may be more complicated, or harder to compute, or may introduce Runge-type oscillations  (\textit{cf.} \cite{Markakis:2014nja}),  we  allow the number of  jumps enforced  to be a specifiable parameter. If we drop jumps in derivatives higher than order $M \in [-1,N]$, then the remaining coefficients are determined by the $N+1$  jump conditions
\begin{equation} \label{fmJ0conditions2}
p^{(m)}(\xi^+)- p^{(m)}(\xi^-)  =   p_ + ^{(m)}(\xi ) - p_ - ^{(m)}(\xi ) = 
    \left\{ \begin{array}{l} {J_m},\quad m = 0,1,...,M\\
    0,\quad m = M + 1,...,N
    \end{array} \right ..
\end{equation}
 Lagrange interpolation is recovered when $M=-1$ or $J_m=0$. 

The ansatz \eqref{PolyN2}  does not amount to domain decomposition as it does not consist of two independent interpolating polynomials of order $n$ and $N-n$ matched at the discontinuity. Our ansatz  covers a single domain $[a,b]$ with  a single piecewise  polynomial of order $N$ which  has piecewise constant coefficients with  a known jump  at the discontinuity $\xi$. 

Substituting Eq.~\eqref{PolyN2} into the collocation conditions \eqref{PolyNConds2} and jump conditions \eqref{fmJ0conditions2} determines the coefficients  $c_j^\pm$. Their explicit form is provided in  Ref.~\cite{Markakis:2014nja}.  However, it is more convenient  to use a  Lagrange basis, whence the piecewise polynomial interpolant reads:   
\begin{equation} \label{PolyNLagrange2}
    p(x) = \sum\limits_{j = 0}^N {{[f_j+\Delta(x_j-\xi;x-\xi) ]}{\pi _j}(x)},
\end{equation}
where
\begin{equation} \label{PolyNLagrangeBasis1}
{\pi _j}(x) = \prod\limits_{\genfrac{}{}{0pt}{2}{k = 0}{
k \ne j}}^N {\frac{{x - {x_k}}}{{{x_j} - {x_k}}}} 
\end{equation}
are the Lagrange basis polynomials. The basis polynomials satisfy
the standard conditions $\pi_j(x_i)=\delta_{ij}$, with $\delta_{ij}$ denoting the Kronecker symbol, so that $p(x)$  satisfies the collocation and jump conditions by construction. The piecewise polynomial interpolant depends  on the location $\xi$ of the discontinuity through the 2-point functions
\begin{equation} \label{sjofx2}
    \Delta(x_j-\xi;x-\xi) = [\theta (x - \xi )\theta (\xi  - {x_j}) - \theta (\xi  - x)\theta ({x_j} - \xi )]\kappa(x_j-\xi) 
\end{equation}
where
\begin{equation} \label{gjofx2}
    \kappa(x_j-\xi) :=  \sum\limits_{m = 0}^{M } {\frac{{{J_m}}}{{m!}}{{({x_j} - \xi )}^m}} 
\end{equation}
are weights computed from the jump conditions at the discontinuity $\xi$ given the nodes $x_j$. As expected, Lagrange interpolation is recovered when no jumps are present, $K_m = 0$. 

 Generally, the nodes and differentiation matrices can be computed and stored in memory. However, in the context of time-dependent partial differential equations with distributional sources, the location and magnitude of the discontinuity will generally change in every time step as the particle moves. Thus, it will be important to be able to evaluate Eqs.~\eqref{sjofx2}
and (\ref{gjofx2}) efficiently. As Eq.~(\ref{gjofx2})  is polynomial in $\delta \xi_j = x_j-\xi$, it can be computed with  higher computational efficiency and numerical precision in Horner form:
\begin{equation} \label{gjofx3}
    \kappa( \delta \xi) =  J_0+\delta \xi \left( J_1+\delta \xi \left(\frac{J_2}{2!}+\delta \xi \left(...+\delta \xi\left(\frac{J_{M-1}}{(M-1)!}+\delta \xi \frac{J_{M}}{M!}\right)...\right)\right)\right ).
\end{equation}
This operation is  vectorized and parallelized across all available cores for all components $\delta \xi_i$ of $\delta \vec \xi$.
When evaluated on a node $x=x_i$, Eq.~\eqref{sjofx2} simplifies to
\begin{equation} \label{sjofx2b}
    \Delta(x_j-\xi;x_i-\xi)  = [\theta (x_{i} - \xi ) - \theta ({x_j} - \xi )]\kappa(x_j-\xi).
\end{equation}
Since the prefactor in $\Delta(x_j-\xi;x_i-\xi) $ is antisymmetric in $i,j$ and $\pi_j(x_i)=\delta_{ij}$ is symmetric, the correction $\sum_j \Delta(x_j-\xi;x-\xi) \pi_j(x) $ in the interpolation formula  \eqref{PolyNLagrange2} vanishes at each node $x=x_i$. This ensures that the collocation conditions \eqref{PolyNConds2} are satisfied.

\subsection{Discontinuous differentiation}  \label{sec:discdif}
Differentiating the piecewise polynomial \eqref{PolyNLagrange2} yields finite-difference or pseudospectral approximations to the $n$-th derivative of $f(x)$. Evaluated at a node $x=x_i$, this yields
\begin{equation} \label{fnderivativeapprox2}
    {f^{(n)}}({x_i}) \simeq {p^{(n)}}({x_i }) = \sum\limits_{j = 0}^N D_{ij}^{(n)}[f_j+ \Delta(x_j-\xi;x_i-\xi)],
\end{equation}
with the differentiation matrices $D_{ij}^{(n)}$ given  by 
\begin{equation}  \label{derivativematricespi}
    D_{ij}^{(n)} = \pi _j^{(n)}({x_i})=\frac{d^n \pi_j(x)}{dx^n} \bigg|_{x=x_i}
\end{equation}
It has been shown \cite{Markakis:2014nja} that at least $M^{\rm th}$-order convergence can be attained when the 2-point $\Delta$ functions  \eqref{sjofx2b} are included in Eq.~\eqref{fnderivativeapprox2}. 

As mentioned above, for time-domain problems, as the discontinuity at $\xi$ moves, one must efficiently update these correction terms in each time-step. Substituting Eq.~\eqref{sjofx2b} into \eqref{fnderivativeapprox2} yields the  expression
\begin{equation} \label{fnderivativeapprox3}
    {f^{(n)}}_i \simeq \sum\limits_{j = 0}^N D_{ij}^{(n)}(f_j+\kappa_j \theta_{i}-\kappa_j \theta_j),
\end{equation}
where $f_i=f(x_i), {f^{(n)}}_i={f^{(n)}}({x_i})$, $\theta_i=\theta (x_{i} - \xi )$ and $\kappa_j=\kappa(x_j-\xi)$ are vectors formed from the values of the respective functions on the set of grid-points $\{x_i\}$. We remark that the above formulae are valid for both pseudospectral and finite-difference methods, \textit{cf.} \cite{Markakis:2014nja} for explicit expressions of the differentiation matrices.

On modern CPUs and GPUs, Eq.~\eqref{fnderivativeapprox3} can be evaluated efficiently via inner (dot product) matrix-vector multiplication and elementwise vector-vector multiplication (Hadamard product). For instance, in Wolfram Language, the simple command
\begin{verbatim}
fn = Dn.f + (Dn.k)*th - Dn.(k*th)
\end{verbatim}
(with \verb|f| and \verb|fn|  denoting $f_i$ and ${f^{(n)}}_i$, \verb|Dn| denoting the $n^{\rm{th}}$ order differentiation matrix given by Eq.~\eqref{derivativematricespi}, and \verb|th|, \verb|k| denoting the vectors  $\theta_i$, $\kappa_i$ respectively) uses the \textsc{Intel Math Kernel Library } to automatically perform the linear algebra operations in Eq.~(\ref{fnderivativeapprox3}) in parallel, across all available cores. Similar commands can be used to accelerate the computation of these products on Nvidia GPUs using cuBLAS.

\section{Discontinuous time-symmetric integration schemes}\label{sec:ch5DiscMethods}

The methods outlined in the previous section enable us to  account for the non-analytic behavior of the target function at a single point along the spatial axis. However, 
the function is also non-analytic when the point is approached along the temporal axis.
This indicates that standard time-steppers (e.g. Runge-Kutta methods) cannot be readily applied to this problem, as they assume the target function to be smooth. Our discontinuous method-of-line experiments   \cite{Markakis:2014nja} indeed showed that using ordinary time-steppers with the (spatially) discontinuous collocation method described above works accurately if the particle is static or stationary  (on a circular orbit), but precision is lost and the method fails to converge when the particle is moving in the $x$ direction. This is in keeping with a key observation made  by Harms et al. \cite{Harms_2014}  that  finite-difference representations of a Dirac-$\delta$ function (like discontinuous collocation methods) converge for a stationary discontinuity but not for a moving discontinuity. 
We will show here that the paradigm of using undetermined coefficients to accommodate known jump and collocation conditions can also be used to develop discontinuous time integration schemes. 
In \S 3.2 below we re-derive the second order discontinuous trapezium rule of \cite{Markakis:2014nja} and in \S 3.3 we extend it to fourth order, obtaining a discontinuous generalization of the Hermite rule.

\subsection{Piecewise smooth time integration}
For generality, we consider the first order differential equation
\begin{equation}
    \frac{dy}{dt} = f(t,y)
\end{equation}
on a small time interval $[t_1,t_2]$. This is equivalent to writing
\begin{equation}\label{eq:ch5IntEq}
    y(t_2) - y(t_1) = \int_{t_1}^{t_2} f(t,y) dt,
\end{equation}
so the problem is now to approximate the above integral. One approach is to construct a polynomial approximant to $f$ between $t_1$ and $t_2$. If the target function $y(t)$ were smooth, it would be straightforward to construct a smooth approximant to $f$ and integrate it to obtain a time-symmetric approximation to the integral, as shown in Ref. \cite{markakis2019timesymmetry}.

Now, though, we shall assume that $y(t)$ is non-analytic at some point $ t_{\rm c}  \in (t_1,t_2)$ and that the discontinuities in $f$ and all its derivatives at this point are known (which has been shown to be the case for linear distributionally-sourced PDEs \cite{Markakis:2014nja}). Let $K_n$ denote the discontinuity in the $n^{\mathrm{th}}$ time derivative of $f$ at $t= t_{\rm c} $ when approached from from below. We accommodate this behavior by constructing the interpolant to be a  \textit{piecewise polynomial} of the form in Eq. \eqref{PolyN2DiscInterpol}, that is,
\begin{equation} \label{eq:piecewisetime}
    f(y,t) \simeq p(t) = \theta(t -  t_{\rm c} ) p_+(t) + \theta( t_{\rm c}  - t) p_-(t).
\end{equation}
Integrating the above interpolant yields discontinuous integration schemes.
We provide second and fourth order symmetric integration schemes below.

\subsection{Discontinuous trapezium rule}
As discussed in Ref.~\cite{Markakis:2014nja}), the simplest interpolant which satisfies endpoint collocation conditions (\textit{e.g.} $f(t_1) = f_1$ and $f(t_2) = f_2$) is a first order polynomial. It follows that the simplest interpolant that satisfies collocation conditions and the jump conditions $K_n$ while preserving time symmetry is one that is first order to both the left and right of $\tau$. Such an interpolant has four parameters so
\begin{equation}
    p_+(t) = A + B t, \quad p_-(t) = C + D t
\end{equation}
Following the construction of the trapezium rule, we impose collocation at the endpoint, so $p_-(t_1) = f_1$ and $p_+(t_2) = f_2$. We satisfy the remaining two degrees of freedom with jump conditions, $p_+( t_{\rm c} ) - p_-( t_{\rm c} ) = K_0$, $p_+'( t_{\rm c} ) - p_-'( t_{\rm c} ) = K_1$. The resulting interpolant is
\begin{subequations}
\begin{align}
    p_+ (t) &= f_1 ~ \frac{t_2 - t}{\Delta t} + f_2 ~ \frac{t - t_1}{\Delta t} + \frac{t_2 - t}{\Delta t}( K_0 - K_1 ( t_{\rm c}  - t_1) )\\
    p_- (t) &= f_1 ~ \frac{t_2 - t}{\Delta t} + f_2 ~ \frac{t - t_1}{\Delta t} + \frac{t - t_1}{\Delta t}( - K_0 - K_1 (t_2 - t_{\rm c} ) )
\end{align}
\end{subequations}
Integrating the piecewise linear polynomial \eqref{eq:piecewisetime} over the time interval $[t_1,t_2]$ yields a discontinuous generalization of the trapezium rule:
\begin{equation}\label{eq:discTrap}
    y_2-y_1 \simeq  \frac{\Delta t}{2} ( f_1 + f_2) + K_0 ~ \frac{\Delta t - 2 \Delta t_{\rm c}}{2} + K_1 ~ \frac{\Delta t_{\rm c}}{2} ( \Delta t_{\rm c} - \Delta t ),
\end{equation}
where $\Delta t_{\rm c} = \tau - t_1$. We shall term this method DH2. The standard trapezium rule is obviously recovered when $K_0 = K_1 = 0$. Although we have not constructed an explicit remainder term, we find that, in practice, this expression exhibits second-order convergence in $\Delta t$.

\subsection{Discontinuous Hermite rule}
The next simplest interpolant that satisfies endpoint collocation is a third-order polynomial. It follows that the next simplest discontinuous piecewise interpolant should be cubic to both the left and right of $\tau$: 
\begin{equation}
    p_+(t) = A + B t + C t^2 + D t^3, \quad p_-(t) = E + F t + G t^2 + H t^3
\end{equation}
Following the same procedure as before, we impose the endpoint collocation conditions $p_-(t_1) = f_1$, $p_+(t_2) = f_2$, osculation conditions $p_-'(t_1) = f_1'$, $p_+'(t_2) = f_2'$, and  jump conditions $p_+( t_{\rm c} ) - p_-( t_{\rm c} ) = K_0$, $p_+'( t_{\rm c} ) - p_-'( t_{\rm c} ) = K_1$, $p_+''( t_{\rm c} ) - p_-''( t_{\rm c} ) = K_2$, and $p_+'''( t_{\rm c} ) - p_-'''( t_{\rm c} ) = K_3$. The resulting interpolant is
\begin{subequations}
\begin{multline}
    p_+(t)=\frac{{{{({t_2}-t)}^2}(2t + {t_2} - 3{t_1})}}{{\Delta {t^{\rm{3}}}}}{f_1} -\frac{{(2t + {t_1} - 3{t_2}){{(t - {t_1})}^2}}}{{\Delta {t^{\rm{3}}}}}{f_2} +\frac{{{{({t_2} - t)}^2}(t - {t_1})}}{{\Delta {t^{\rm{2}}}}}{{f}'_1}- \\
    \frac{{({t_2} - t){{(t - {t_1})}^2}}}{{\Delta {t^{\rm{2}}}}}{f'^{}_2} +\frac{(t_2 - t)^2 (2 t-3 t_1+t_2)}{\Delta t^3} K_0\\
    - \frac{(t-t_2)^2 (\Delta t_{\rm c}  (2 t-3 t_1 + t_2) - (t-t_1) \Delta t )}{\Delta t^3} K_1\\
    + \frac{(t_2-t)^2 \Delta t_{\rm c} (\Delta t_{\rm c}  (2 t-3 t_1 +t_2 )- 2 (t-t_1) \Delta t)}{2 \Delta t^3} K_2\\ - \frac{\Delta t_{\rm c}^2 (t_2 - t)^2 (\Delta t_{\rm c}  (2 t-3 t_1 + t_2)-3 (t-t_1) \Delta t)}{6 \Delta t^3} K_3
\end{multline}
\begin{multline}
    p_-(t)=\frac{{{{({t_2}-t)}^2}(2t + {t_2} - 3{t_1})}}{{\Delta {t^{\rm{3}}}}}{f_1} -\frac{{(2t + {t_1} - 3{t_2}){{(t - {t_1})}^2}}}{{\Delta {t^{\rm{3}}}}}{f_2} +\frac{{{{({t_2} - t)}^2}(t - {t_1})}}{{\Delta {t^{\rm{2}}}}}{{f}'_1}\\
    - \frac{{({t_2} - t){{(t - {t_1})}^2}}}{{\Delta {t^{\rm{2}}}}}{f'^{}_2} +\frac{(t - t_1)^2 (2 t-3 t_2+t_1)}{\Delta t^3} K_0\\
    - \frac{(t - t_1)^2 (\Delta t_{\rm c}  (2 t-3 t_2 + t_1) + (2t_2 - t - t_1) \Delta t )}{\Delta t^3} K_1\\
    + \frac{(t-t_1)^2 (t_2 - \tau) (\Delta t_{\rm c}  (2 t-3 t_2 +t_1 )- 2 \Delta t^2)}{2 \Delta t^3} K_2\\
    + \frac{(t - t_1)^2 (t_2 - \tau)^2 (\Delta t_{\rm c}  (2 t-3 t_2 + t_1)- (t-t_1) \Delta t)}{6 \Delta t^3} K_3
\end{multline}
\end{subequations}
Integrating the   piecewise cubic polynomial \eqref{eq:piecewisetime} over the interval $[t_1,t_2]$ yields a discontinuous generalization of the Hermite rule:
\begin{multline}\label{eq:discHerm}
    y_2 - y_1 \simeq \frac{\Delta t}{2}(f_1 + f_2) + \frac{\Delta t^2}{12}(f_1' - f_2') + K_0 ~ \frac{\Delta t - 2 \Delta t_{\rm c}}{2}  + K_1~ \frac{\Delta t^2 - 6 \Delta t \Delta t_{\rm c} + 6 \Delta t_{\rm c}^2}{12} \\
    - K_2 ~ \frac{\Delta t_{\rm c}(\Delta t^2 - 3 \Delta t \Delta t_{\rm c} + 2 \Delta t_{\rm c}^2)}{12} + K_3~ \frac{\Delta t_{\rm c}^2 ( \Delta t - \Delta t_{\rm c})^2}{24} 
\end{multline}
We shall term this method DH4. The standard Hermite rule is recovered when the $K_n$ are set to zero. As before, we find that, in practice, this expression exhibits fourth order convergence in $\Delta t$ even though we have not computed a remainder.

\subsection{Numerical experiment}
We demonstrate the validity of these discontinuous integration schemes by studying how well they approximate the definite integral of a discontinuous function:
\begin{equation}\label{eq:ch5ToyFunc}
    f(t) = P_5(t) \theta(t) + Q_3(t) \theta(-t)
\end{equation}
where $P_5$ is the fifth Legendre polynomial and $Q_3$ is the third associated Legendre polynomial. The first four jumps, obviously located at $t=0$, may be calculated analytically:
\begin{equation}
    K_0 = -\frac{2}{3}, \quad K_1 = \frac{15}{8}, \quad K_2 = 8, \quad K_3 = - \frac{105}{2}
\end{equation}
We approximate the integral of this function over the interval $[-0.51,0.49]$ (this interval was chosen so that $t=0$ is never located on a grid point, which Eqs. \eqref{eq:discTrap} and \eqref{eq:discHerm} implicitly assume), which may be calculated analytically to be $\approx 0.294859$.

\begin{figure}[!t]
    \centering
    \begin{subfigure}{0.45\textwidth}
        \centering
        \includegraphics[width=\textwidth]{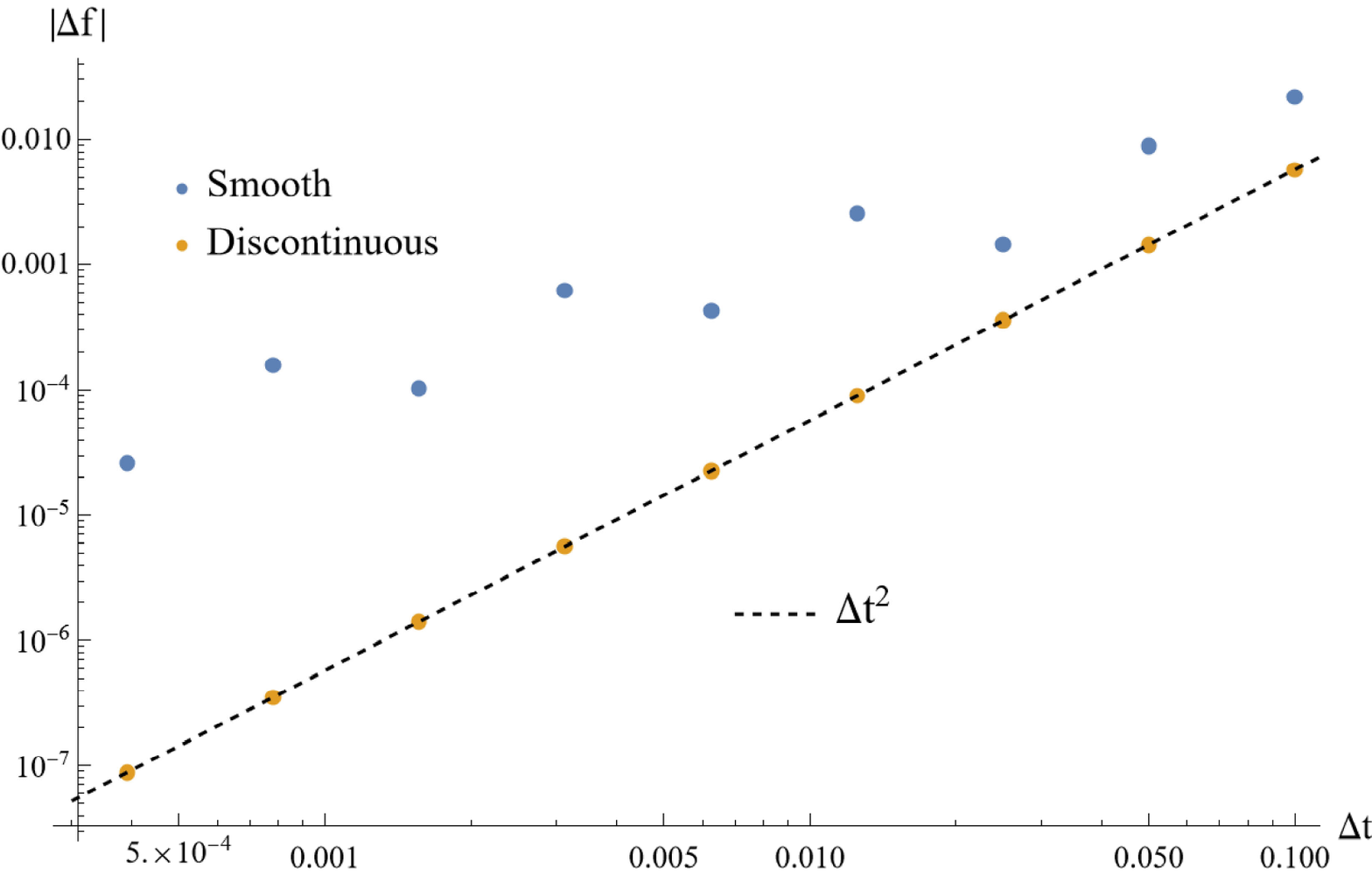}
        \caption{Trapezium}
    \end{subfigure}
    \hfill
    \begin{subfigure}{0.45\textwidth}
        \centering
        \includegraphics[width=\textwidth]{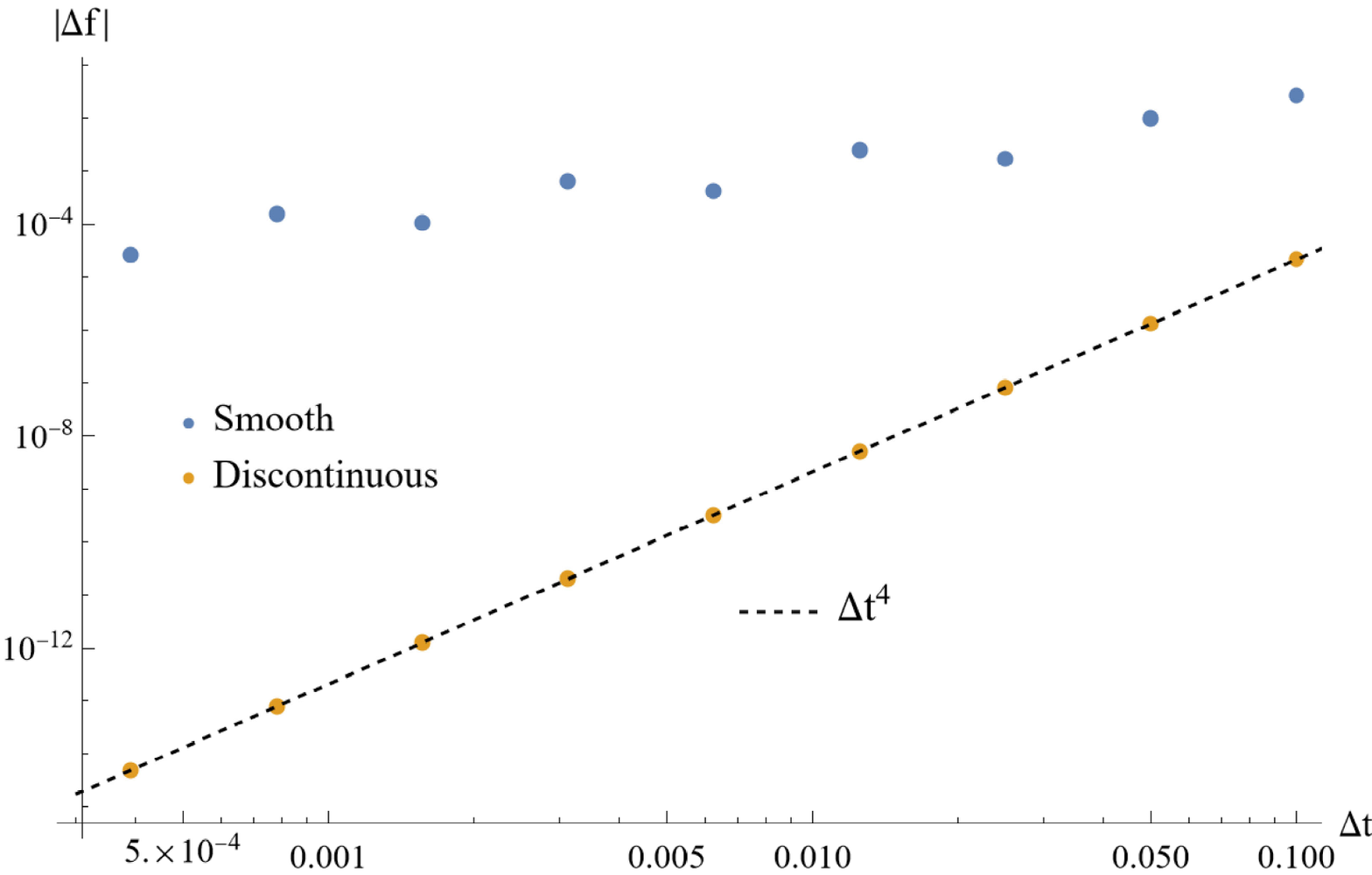}
        \caption{Hermite}
    \end{subfigure}
    \caption{The error resulting from using Eqs.~\eqref{eq:discTrap} (on the left) and \eqref{eq:discHerm} (on the right) to approximate the integral of Eq.~\eqref{eq:ch5ToyFunc} on the interval $(-0.51,0.49)$. We observe that Eq. \eqref{eq:discTrap} exhibits second order convergence while Eq. \eqref{eq:discHerm} exhibits fourth order convergence while the smooth versions of these schemes result in poorly behaved approximations.}
    \label{fig:ch5DiscError}
\end{figure}

In Figure \ref{fig:ch5DiscError}, We observe that the approximation resultingfrom applying Eq. \eqref{eq:discTrap} to the subinterval containing $t=0$ indeed exhibits error $\mathcal{O}(\Delta t^2)$ and Eq.~\eqref{eq:discHerm} exhibits error $\mathcal{O}(\Delta t^4)$, as we might expect from their construction. In addition, we note that the smooth versions of these expressions, without accounting for the jumps, result in poorly-behaved approximations.

\section{Application to distributionally-sourced PDEs}\label{sec:ch5deltaPDEs}

\subsection{Scalar wave equation with distributional source term}\label{sec:FlatDeltaPDE}
To illustrate the utility of both the discontinuous time stepping and spatial collocation schemes, we consider a simple prototype of a distributionally-sourced wave equation, and a variant of Eq.~\eqref{eq:DistrForcedPDE},  considered previously in \cite{Field2009,Markakis:2014nja}: 
\begin{equation}\label{eq:FieldFlat}
    - \partial^2_t \Psi + \partial^2_x \Psi = F(t) ~ \delta'( x - \xi( t )) +G(t) ~ \delta( x - \xi( t )).
\end{equation}
Heuristically, this equation gives the field produced by a scalar charge  moving along the worldline $\xi(t)$ with a time-dependent monopole and dipole moment.
 
 Eq.~\eqref{eq:FieldFlat} possesses closed form solutions  \cite{Field2009,Markakis:2014nja} if the particle is taken to move at a constant speed $v$ along a worldline $\xi(t) = v \,t$. If $F(t) = 0$ and $G(t) = \cos t$, ``Solution I'' is given by
\begin{equation}\label{eq:ch5FieldI}
    \Psi_{\rm I}(t,x) = - \frac{1}{2} \sin \Big( \gamma^2(t-vx - |x - v t|) \Big)
\end{equation}
If $F(t) = \sin t$ and $G(t) = 0$, ``Solution II'' is given by
\begin{equation}\label{eq:ch5FieldII}
    \Psi_{\mathrm{II}}(t,x) = \frac{1}{2} \gamma^2(v + \mathrm{sgn}(x - v t)) \cos \Big( \gamma^2(t-vx - |x - v t|) \Big)
\end{equation}
In prior work  \cite{Markakis:2014nja}, the above solutions were used to test numerical solutions to Eq.~\eqref{eq:FieldFlat}  in  Minkowski coordinates $\{t,x\}$ using discontinuous time symmetric and discontinuous collocation methods. This required imposing boundary conditions at the ends of the spatial domain, which is generally not straightforward when a potential term is added to the above equation (as is necessary, for instance, in black hole perturbation theory). Here to automatically impose boundary conditions and improve computational efficiency, we adopt these numerical methods to hyperboloidal coordinates.

\subsection{Hyperboloidal slicing}

We adopt the method  of \textit{hyperboloidal compactification} developed by Zenginoglu \cite{zenginoglu_hyperboloidal_2008,PhysRevD.83.127502,Zengino_lu_2009,ZENGINOGLU20112286}. As demonstrated in \cite{moboylethesis,Markakis:2023pfh}, the idea is as follows. In Minkowski space, the null rays (characteristic curves of the scalar wave equation) form lines that are $45^\circ$ from the $t$- and $x$-axes. If instead one defines a coordinate by
\begin{equation}
    \tau = t - h(x)
\end{equation}
and demands that
\begin{equation}
    |h'(x)|<1, \quad -\infty < x < \infty,
\end{equation}
the new coordinate will be timelike throughout the interior of the domain. Moreover, if $h$ satisfies
\begin{equation}
    \lim_{x\to \pm \infty} h'(x) = \mp 1,
\end{equation}
then the coordinate $\tau$ will become null on the boundary surfaces and therefore intersect   $\mathcal{I}^+$. Such a coordinate is termed \textit{hyperboloidal}. Once this is achieved, one can bring the boundary surfaces into a finite domain via compactification; i.e. choose a new spacelike coordinate
\begin{equation}
    x=g(\sigma) 
\end{equation}
such that
\begin{equation}
    \sigma( [- \infty, \infty]) = [0, 1]
\end{equation}
Thus, such a set of coordinates maps the behaviors of the field on  $\mathcal{I}^+$ to finite spacelike and interior-timelike coordinates and preserves the initial value formulation for initial data given on a constant $\tau$ slice $\Sigma_\tau$, making it ideal for numerical studies of GW generation.

Several such coordinate slicings have been found for the Minkowski, Schwarzschild and Kerr  spacetime \cite{zenginoglu_hyperboloidal_2008,racz_numerical_2011,Gautam:2021ilg}, but we find that the ``minimal gauge'' defined by Ansorg and Macedo \cite{ansorg_spectral_2016} yields the simplest algebraic expressions and covers the entire black hole exterior with a single hyperboloidal layer. As written, the following compactification can be used in Minkowski as well as Schwarzschild spacetime (with $t,x$  taken to be tortoise coordinates):
\begin{eqnarray}\label{eq:compact_g}
g(\sigma) \mkern-10mu &:=& \mkern-10mu \int \frac{{1 }}{{2{\sigma ^2}( \sigma -1)}}
= \frac{1}{2} \left(\frac{1}{\sigma }+\log (1-\sigma )-\log (\sigma )\right) \quad  \quad
\end{eqnarray}
In \cite{Markakis:2023pfh}, the asymptotic behavior of null rays was used to show that a suitable height function choice is:
\begin{equation}\label{eq:height_h}
h(\sigma) := g(\sigma) - \frac{1}{\sigma } + \ln \sigma +\mathcal{O}(\sigma).  \end{equation}
Truncating this expression to linear order in $\sigma$ amounts to the so-called ``minimal gauge''. These expressions can be further generalized to Kerr spacetime \cite{Markakis:2023pfh}.

Then, \eqref{eq:FieldFlat} can be written  in the form
\begin{multline}\label{eq:ch6GeneralWave}
    Z(\sigma) \partial_\tau^2 \Psi + A(\sigma) \partial_\tau \partial_\sigma \Psi + B(\sigma) \partial_\tau \Psi + C(\sigma) \partial_\sigma^2 \Psi  + E(\sigma) \partial_\sigma \Psi \\
    = F(\tau,\sigma) \delta'(\sigma - \xi(\tau)) + G(\tau,\sigma) \delta(\sigma-\xi(\tau))
\end{multline}
Eq.~\eqref{eq:ch6GeneralWave} can be obtained from the 1+1 covariant wave equation,
\begin{equation}\label{eq:LittlePhiEqnCovar}
\Box_{\mathcal{N}^2} \Psi :=\eta^{\alpha \beta} \nabla_{\alpha} \nabla_{\beta} \Psi   = S,
\end{equation}
on a flat manifold $\mathcal{N}^2$ spanned by $\tau$ and $\sigma$, 
where $\eta_{\alpha \beta}$ is an effective Minkowski metric \cite{jaramillo_pseudospectrum_2021} given by
 the line element,
\begin{equation} \label{eq:hyperbMink}
    %ds^2 =
    \eta_{\alpha \beta}dx^\alpha dx^\beta = - d\tau^2 + \frac{1 - 2 \sigma^2}{\sigma^2 (1-\sigma)} ~ d\tau d \sigma + \frac{1+\sigma}{\sigma^2(1-\sigma)} d\sigma^2.
\end{equation}
In the above coordinate chart, the wave operator for a scalar field  is polynomial in $\sigma$ and  given by
\begin{multline}\label{eq:ch6FlatHyperEqn}
    \Box_{\mathcal{N}^2} \Psi = - 4 \sigma^2(1-\sigma^2) \partial_\tau^2 \Psi + 4\sigma^2(1-\sigma)(1-2\sigma^2) \partial_\tau \partial_\sigma \Psi + 4 \sigma^4 (1-\sigma)^2 \partial_\sigma^2 \Psi \\
    - 8 \sigma^3(1-\sigma) \partial_\tau \Psi + 4\sigma^3 (1-\sigma)(2-3\sigma) \partial_\sigma \Psi,
\end{multline}
which allows the coefficients $A,B,C, E$ and $Z$ to be read off:
\begin{equation}\label{eq:ScalarL}
    A(\sigma) = {1 - 2 \sigma^2}, \quad B(\sigma) = - {2 \sigma}, \quad 
    C(\sigma) = {1 - \sigma}, \quad E(\sigma) = {\sigma(2-3\sigma)}, \quad Z(\sigma) = {1+\sigma}.
\end{equation}

To make manifest that outflow boundary conditions are automatically imposed at the  boundaries $\sigma=0$ and $\sigma=1$ in this coordinate chart, we calculate the characteristic speeds of this wave equation,
\begin{equation}
    v_+ (\sigma)= \frac{\sigma^2}{1+\sigma}, \quad v_- (\sigma)= -1 + \sigma.
\end{equation}
As required, the light cones tip over at the two boundaries, and we have $v_+ (0)=0$ and $v_- (1)=0$.

 In the slicing of Eq.~\eqref{eq:hyperbMink}, the scalar solutions \eqref{eq:ch5FieldI}-\eqref{eq:ch5FieldII} of Eq.~\eqref{eq:FieldFlat} transform to read
\begin{equation}\label{eq:ch6FieldIHyper}
    \Psi_{\rm I}(\tau,\sigma) = 
    \begin{cases}
        - \frac{1}{2} \sin \bigg( \frac{1 + \sigma (\tau - \log \sigma)}{\sigma (1+v)} \bigg) & \sigma > \zeta(\tau)\\
        - \frac{1}{2} \sin \bigg( \frac{ \tau - \log(1 - \sigma)}{1-v} \bigg) & \sigma < \zeta(\tau)
    \end{cases},
\end{equation}
\begin{equation}\label{eq:ch6FieldIIHyper}
    \Psi_{\rm II}(\tau,\sigma) = 
    \begin{cases}
        - \frac{1}{2(1+v)} \cos \bigg( \frac{1 + \sigma (\tau - \log \sigma)}{\sigma (1+v)} \bigg) & \sigma > \zeta(\tau)\\
        \frac{1}{2(1-v)} \cos \bigg( \frac{ \tau - \log(1 - \sigma)}{1-v} \bigg) & \sigma < \zeta(\tau)
    \end{cases},
\end{equation}
where %$\xi(\tau)=g(\zeta(t))$ 
%is
$\zeta(\tau)=g^{-1}(\xi(t))$ 
is the worldline of the particle in the hyperboloidal chart $\{\tau,\sigma \}$. This provides exact solutions by which to study the converge of our discontinuous time steppers below.

\begin{figure}[!t]
    \centering
     \includegraphics[width=\textwidth]{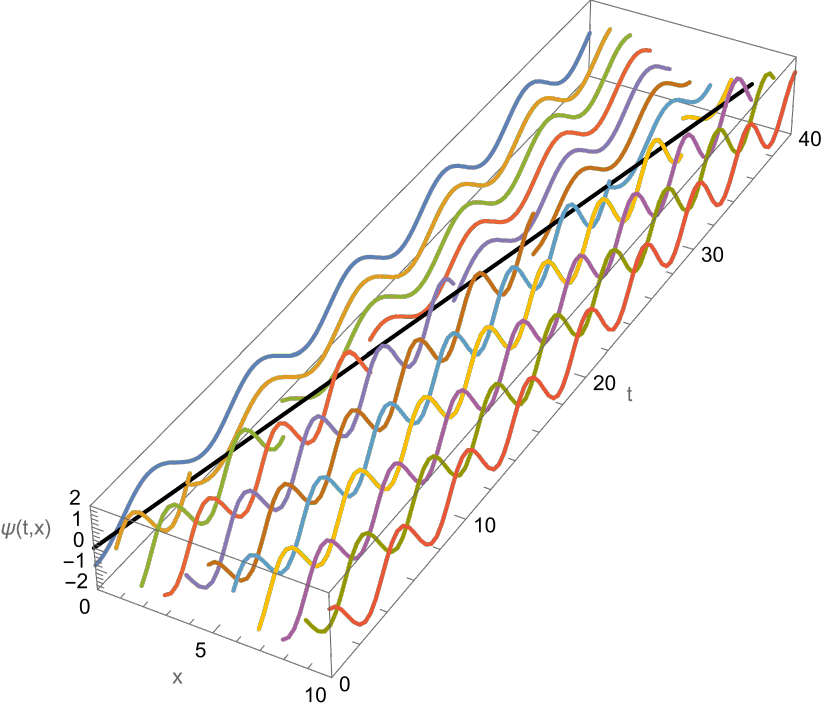}
    \caption{The solution~\eqref{eq:ch5FieldII} (or, equivalently, Eq.~\eqref{eq:ch6FieldIIHyper} in hyperboloidal coordinates), recovered with a discontinuous method of lines. The black line represents the wordline $\xi(t) =v \, t$ of the particle. Each colored line represents the solution $\Psi(t,x)$ on the respective grid-point $x=x_i$ as a function of time $t$. The solution $\Psi(t,x)$ (and its spatial and temporal derivatives) is discontinuous accross the worldline $\xi(t)$, both in the space and time direction. Thus, a discontinuous collocation method is required to differentiate $\Psi$ in the $x$ direction, and a discontinuous time integration method is required to integrate the solution in the $t$ direction. Specifically, every time the particle worldline crosses a grid-point, the time integrator must account for the jump in the solution and its time derivatives, as described in Sec.~3.}
    \label{fig:DiscMOL}
\end{figure}

\subsection{Generalized recursion relation}
Before we turn the numerical evolution of specific problems, we first address the question of how jumps in the target function may be calculated in arbitrary coordinate systems. The method of unit jump functions presented in \cite{Markakis:2014nja} may be extended to this situation, but it is significantly more involved than before.

Eq.~\eqref{eq:FieldFlat} has been studied by Field et al.  \cite{Field2009} and Markakis et al. \cite{Markakis:2014nja}. Here, to develop the computational tools necessary for future hyperboloidal black hole perturbation applications, we study this equation on a hyperboloidal slice~\eqref{eq:hyperbMink}, given by Eq.~\eqref{eq:ch6GeneralWave}. As in earlier work \cite{Markakis:2014nja},  we suppose that the general solution to Eq.~\eqref{eq:ch6GeneralWave} may be written as the superposition $\Psi = \Psi_{\rm H} + \Psi_{\rm NS}$, where $\Psi_{\rm H}$ is a general solution to the homogeneous equation and $\Psi_{\rm NS}$ is a particular solution corresponding to the source terms. We note that, since the source terms are distributions, the part $\Psi_{\rm NS}$ of the target function is necessarily non-smooth. We begin with the same decomposition of the non-smooth part of the target function,
\begin{equation}\label{eq:ch6NSFunc}
    \Psi_{\rm NS}(\tau,\sigma) = \sum_{n=0}^\infty J_n(\tau) \Phi_n(\sigma;\zeta),
\end{equation}
but now assume the more general form \eqref{eq:ch6GeneralWave} for the evolution equation.

We substitute the form Eq.~\eqref{eq:ch6NSFunc} into this problem and simplify the terms of the left hand side of Eq.~\eqref{eq:ch6GeneralWave}. In evaluating the source terms on the right hand side, we now consider the possibility of $F$ and $G$ depending on $\sigma$ as well as $\tau$, so we must invoke the $\delta$ function \textit{selection properties} (\textit{cf.} Appendix D of \cite{moboylethesis}):
\begin{eqnarray}
\label{eq:deltaPrimeSelection}
    &&f(\sigma) \delta'(\sigma-\sigma_0) = f(\sigma_0) \delta'(\sigma-\sigma_0) - f'(\sigma_0) \delta(\sigma-\sigma_0)\\
\label{eq:deltaSelection}
    &&g(\sigma) \delta(\sigma- \sigma_0) = g(\sigma_0) \delta(\sigma- \sigma_0)
\end{eqnarray}
Applied to this problem, we have
\begin{multline}
    F(\tau,\sigma) \delta'(\sigma - \zeta(\tau)) + G(\tau,\sigma) \delta(\sigma-\zeta(\tau))\\
    = F(\tau,\zeta(\tau)) \delta'(\sigma - \zeta(\tau)) + \Big( - \partial_\sigma F(\tau,\zeta(\tau)) + G(\tau,\zeta(\tau)) \Big) \delta(\sigma - \zeta(\tau))
\end{multline}
The details of this calculation are presented in Appendix \ref{ch:recursion}. The end results are explicit expressions for the first two jumps,
\begin{equation}\label{eq:ch6NewJ0}
    J_0 = \gamma^2 F(\tau,\zeta)
\end{equation}
\begin{multline}\label{eq:ch6NewJ1}
    J_1 = \gamma^2 \bigg( - \partial_\sigma F(\tau,\zeta) + G(\tau,\zeta) + \Dot{J}_0 \Big( 2 Z(\zeta) \Dot{\xi} - A(\zeta) \Big) \\
    + J_0 \Big( Z'(\zeta) \Dot{\zeta}^2 - A'(\zeta) \Dot{\zeta} + C'(\zeta) + Z(\zeta) \Ddot{\zeta} + B(\zeta) \Dot{\zeta} - E(\zeta) \Big) \bigg)
\end{multline}
where $\gamma$ is the Lorentz factor given by
\begin{equation}
    \gamma^{-2} = Z(\zeta) \Dot{\zeta}^2 - A(\zeta) \Dot{\zeta} + C(\zeta).
\end{equation}
A recursion relation  specifies all higher order jumps:
\begin{multline}\label{eq:ch6NewRecur}
    J_{n+2} = - \gamma^2 \Bigg( \sum_{k=0}^n \binom{n}{k} \Big( A^{(k)}(\zeta) \Dot{J}_{n+1-k} + B^{(k)}(\zeta) (\Dot{J}_{n-k} - \Dot{\zeta} J_{n+1-k}) - E^{(k)}(\zeta) J_{n+1-k} \\
    + Z^{(k)}(\zeta) ( \Ddot{J}_{n-k} - 2 \Dot{J}_{n+1-k} \Dot{\zeta} + J_{n+1-k} \Ddot{\zeta} ) \Big) + \sum_{k=1}^n \binom{n}{k} J_{n+2-k} \Big( A^{(k)}(\zeta) \Dot{\zeta} + C^{(k)}(\zeta) \Big) \Bigg)
\end{multline}

The coefficients $J_n(\tau)$ give the discontinuities in the $n^\mathrm{th}$ derivative of $\Psi_\mathrm{NS}$. In fact, noting that $\Phi_0(\tau,\sigma) = \frac{1}{2} \mathrm{sgn}(\sigma-\zeta(\tau))$,
\begin{equation} \label{eq:space_jump}
    \llbracket \Psi_\mathrm{NS} \rrbracket = \lim_{\sigma \rightarrow \zeta^+} \Psi_\mathrm{NS}(\tau,\sigma) - \lim_{\sigma \rightarrow \zeta^-} \Psi_\mathrm{NS}(\tau,\sigma) = J_0(\tau)
\end{equation}
If instead we take the limits in the $\tau$-direction to obtain the first jump $K_0$, we have
\begin{equation}\label{eq:time_jump}
    K_0 = \lim_{\tau \rightarrow T^+} \Psi_\mathrm{NS}(\tau,\sigma) - \lim_{\tau \rightarrow T^-} \Psi_\mathrm{NS}(\tau,\sigma) = - J_0(T) = - \llbracket \Psi_\mathrm{NS} \rrbracket,
\end{equation}
where we have defined $T = \zeta^{-1}(\sigma)$. Thus, a jump along the $\tau-$direction is numerically the negative of a jump along the $\sigma-$direction. This observation is crucial to the application of discontinuous Hermite integrators to time evolution: once the jumps are calculated in the spatial direction, they are easily adapted to the time direction for use in Eqs.~\eqref{eq:discTrap} and \eqref{eq:discHerm}.

\subsection{Distributionally-sourced wave equation in hyperboloidal coordinates}
We wish to solve Eq.~\eqref{eq:FieldFlat} by employing hyperboloidal compactification; that is, in a coordinate system where $\mathcal{I}^+$, the only null boundary surface in Minkowski space, is brought to finite coordinate values. Our studies of effective action in the previous section has  yielded such a coordinate transformation on Minkowski space, Eq.~\eqref{eq:LittlePhiEqnCovar}. 

The only remaining question is how to handle the source terms in these coordinates. We wish to study Solutions I and II in \cite{Field2009}, so we  consider the particle worldline $\xi(t) = v \; t$ giving rise to the  source terms $F(t) \delta'(x - vt)$ and $G(t) \delta(x-vt)$ in a standard Minkowski coordinate chart $\{t,x\}$. Since the $\delta$-functions are given in another chart, we bring them to the new $\{\tau,\sigma\}$ chart by invoking the $\delta$-function \textit{composition rules} (\textit{cf.} Appendix D of \cite{moboylethesis}): 
\begin{eqnarray}
    \label{eq:ch6DeltaComp}
    &&\delta(f(\sigma)) = \frac{1}{|f'(\sigma_0)|} ~\delta(\sigma-\sigma_0)\\
    \label{eq:ch6DeltaPComp}
    &&\delta'(f(\sigma))= \frac{f'(\sigma_0)}{|f'(\sigma_0)|^3}~ \delta'(\sigma- \sigma_0) + \frac{f''(\sigma_0)}{|f'(\sigma_0)|^3}~ \delta(\sigma- \sigma_0)
\end{eqnarray}
where $f(\sigma_0) = 0$. Noting that $x = x(\sigma)$ and $t=t(\tau,\sigma)$, we note that the $\delta$-function source terms for this problem transform as
\begin{equation}\label{eq:ch6NewDelta}
    G(t)\delta(x - v t) = \frac{G(t_\zeta)}{|\Delta(\zeta)|} \delta(\sigma - \zeta)
\end{equation}
\begin{multline}\label{eq:ch6NewDeltaP}
    F(t) \delta'(x - v t) = \frac{\Delta(\zeta)}{|\Delta(\zeta)|^3} F(t_\zeta) \delta'(\sigma - \zeta)
    + \bigg( \frac{\Delta'(\zeta)}{|\Delta(\zeta)|^3} F(t_\zeta) - \frac{\Delta(\zeta)}{|\Delta(\zeta)|^3} F'(t_\zeta) \frac{\partial t}{\partial \sigma}\Big|_\zeta \bigg)\delta(\sigma- \zeta)
\end{multline}
where we have invoked both the selection and composition properties of the $\delta$-function and we have let $\sigma=\zeta(\tau)$ be the worldline of the particle ($x=v \;t=\xi(t)$) transformed to the new coordinate chart $\{ \tau, \sigma\}$. For convenience we have also defined
\begin{equation}
    \Delta(\zeta) = \frac{\partial}{\partial \sigma}( x - v t) \Big|_{\sigma=\zeta}
\end{equation}
and
\begin{equation}
    t_\zeta = t(\tau,\zeta(\tau)).
\end{equation}
The particle position in the chart \eqref{eq:hyperbMink} is governed by the equation of motion
\begin{equation}
    \frac{d \zeta}{d \tau} = \frac{2 v(1-\zeta)\zeta^2}{1 + v(1-2\zeta^2)}.
\end{equation}
which may be solved to give $\zeta(\tau)$ either implicitly or numerically. We further note that, in this chart,
\begin{equation}
    \Delta(\zeta)  = \frac{1 + v(1-2\zeta^2)}{2\zeta^2 (1-\zeta)}.
\end{equation}
and
\begin{equation}
    t_\zeta = \frac{1 + 2\zeta \tanh^{-1} (1-2\zeta)}{2 v\zeta}
\end{equation}
This information now allows us to completely determine the jumps in both Solutions I (Eq.~\eqref{eq:ch5FieldI} and II (Eq.~\eqref{eq:ch5FieldII}) using the generalized recursion relation. The first few are tabulated in Appendix \ref{ch:Jumps}.

Next, we incorporate the transformed $\delta$-function source terms in the same manner as in Ref. \cite{Markakis:2014nja}. For generality, we suppose that the wave equation is written in the form of Eq.~\eqref{eq:ch6GeneralWave} which, upon first order in time reduction, takes the form:
\begin{equation}  \label{eq:DiscFirstOrder}
\left\{ \begin{gathered}
  {\partial _\tau }\Psi  - \Pi  = \tilde F(\tau )\delta (\sigma  - \zeta) \hfill \\
  Z(\sigma ){\partial _\tau }\Pi  + A(\sigma ){\partial _\sigma }\Pi  + B(\sigma )\Pi  + C(\sigma )\partial _\sigma ^2\Psi  + E(\sigma ){\partial _\sigma }\Psi  = \tilde G(\tau )\delta (\sigma  - \zeta ) \hfill \\ 
\end{gathered}  \right.
\end{equation}
Here, to avoid the difficulty of a $\delta'$ function in a first-order equation, we have defined the variable $\Pi$ by
\begin{equation}
    \partial_\tau \Psi -\Pi = \frac{\Delta(\xi)}{|\Delta( \zeta)|^3} \frac{F(t_ \zeta)}{A( \zeta) - \Dot{ \zeta}Z( \zeta)} \delta (\sigma - \zeta) :=  \Tilde{F}(\tau) \delta (\sigma -  \zeta).
\end{equation}
Substitution of $\Pi$ into the field equation~\eqref{eq:ch6GeneralWave} yields an equation that is first-order in $\Psi, \Pi$ given by
\begin{multline}\label{eq:ch6FullFlatEqn}
    Z(\sigma) \partial_\tau \Pi + A(\sigma) \partial_\sigma \Pi + B(\sigma) \Pi + C(\sigma) \partial_\sigma^2 \Psi + E(\sigma) \partial_\sigma \Psi\\
    = \Bigg( \frac{\Delta'( \zeta)}{|\Delta( \zeta)|^3} F(t _\zeta) - \frac{\Delta( \zeta)}{|\Delta( \zeta)|^3} F'(t_ \zeta) \frac{\partial t}{\partial \sigma}\Big|_ \zeta
    + \frac{G(t_ \zeta)}{|\Delta( \zeta)|}- Z( \zeta) \Tilde{F}'(\tau)\\
    - (A'( \zeta) - \Dot{ \zeta} Z'( \zeta)) \Tilde{F}(\tau) \Bigg) \delta(\sigma -  \zeta)
    := \Tilde{G}(\tau) \delta (\sigma- \zeta)
\end{multline}

We proceed to demonstrate that the discontinuous collocation methods of Section \ref{sec:DiscColl} with the discontinuous time integration methods of Section \ref{sec:ch5DiscMethods} may be used in concert to solve Eq.~\eqref{eq:ch6FlatHyperEqn} with distributional source terms Eqs.~\eqref{eq:ch6NewDelta} and \eqref{eq:ch6NewDeltaP}. Before doing so, we first take the time integration method used in Section \ref{sec:ch5deltaPDEs} and reformulate it using a first-order in time ``state vector'' approach so as to simplify the matrix equations used in evolution. %We first observe that Eq.~\eqref{eq:ch6FlatHyperEqn} may be readily cast in the form of Eqs. \eqref{eq:ch6StateVector} and \eqref{eq:FirstOrderReduced} by defining a regular function $\psi = \Psi / \sigma$, setting $V_\ell=0$ in \eqref{eq:SecondOrderOps}, and retaining the definitions of $A,B,C,$ and $E$ in Eqs. \eqref{eq:ScalarL} and \eqref{eq:ScalarM}. 
Since the target function is now discontinuous, there is the added complication of correcting the spatial derivative operators in Eq.~\eqref{eq:DiscFirstOrder}.

We observe that discontinuous corrections only enter as a consequence of the $A,C$, and $E$ coefficients. Using Eq. \eqref{fnderivativeapprox2}, we note that the above-mentioned terms discretize to
\begin{eqnarray}
    &&E(\sigma)\partial_\sigma \Psi |_{\sigma=\sigma_i} \simeq \sum\limits_{j = 0}^N \Big(\mathbf{E} \; \mathbf{D}^{(1)} \Big)_{ij}[\Psi_j+\Delta(\sigma_j- \zeta;\sigma_i- \zeta)], \quad \mathbf{E} = \mathrm{diag}(E(\sigma_i))\\
    &&C(\sigma)\partial^2_\sigma \Psi |_{\sigma=\sigma_i} \simeq \sum\limits_{j = 0}^N \Big(\mathbf{C}  \; \mathbf{D}^{(2)} \Big)_{ij}[\Psi_j+\Delta(\sigma_j- \zeta;\sigma_i- \zeta)], \quad \mathbf{C} = \mathrm{diag}(C(\sigma_i))\\
    &&A(\sigma)\partial_\sigma \Pi |_{\sigma=\sigma_i} \simeq \sum\limits_{j = 0}^N \Big(\mathbf{A}  \; \mathbf{D}^{(1)} \Big)_{ij}[\Pi_j+\Delta_\Pi(\sigma_j- \zeta;\sigma_i- \zeta)], \quad \mathbf{A} = \mathrm{diag}(A(\sigma_i))
\end{eqnarray}
where an inner (or dot) product between $\mathbf{E}$, $\mathbf{C}$, $\mathbf{A}$ and the differentiation matrices $\mathbf{D}^{(n)}$ is implied.

This indicates that the discretized Eq.~\eqref{eq:DiscFirstOrder} should now include an effective source term of the form
\begin{equation}
    \mathbf{s}(\tau) =
    \begin{pmatrix}
        \mathbf{0}\\
        \mathbf{r}^{A}(\tau) + \mathbf{s}^{C}(\tau) + \mathbf{s}^{E}(\tau)
    \end{pmatrix}
\end{equation}
where
\begin{eqnarray}
    && s^E_i(\tau)= \sum\limits_{j = 0}^N \Big(\mathbf{E}  \; \mathbf{D}^{(1)} \Big)_{ij} \Delta(\sigma_j- \zeta;\sigma_i- \zeta)\\
    && s^C_i(\tau)= \sum\limits_{j = 0}^N \Big(\mathbf{C}  \; \mathbf{D}^{(2)} \Big)_{ij} \Delta(\sigma_j- \zeta;\sigma_i- \zeta)\\
    &&r^A_i(\tau)= \sum\limits_{j = 0}^N \Big(\mathbf{A}  \; \mathbf{D}^{(1)} \Big)_{ij}\Delta_{\Pi}(\sigma_j- \zeta;\sigma_i- \zeta)
\end{eqnarray}
Note that $\Delta_\Pi$ is defined in the same manner as Eq.~\eqref{sjofx2}, but now the $\kappa$ function is calculated using the jumps in $\Pi$ instead of $\Psi$.

Putting everything together, the evolution equation~\eqref{eq:DiscFirstOrder} for the discretized state vector becomes
\begin{equation}\label{eq:ch6GeneralMatrixEqn}
    \frac{d \mathbf{u}}{d \tau} = \mathbf{L} \,\mathbf{u} + \mathbf{s}(\tau) + \mathbf{\Tilde{F}}(\tau) \delta \big( \bm{\sigma} - \zeta \big)
\end{equation}
where
\begin{equation}
 \mathbf{{u}}(\tau) =
    \begin{pmatrix}
       \boldsymbol{{\Psi}}(\tau)\\
        \boldsymbol{\Pi}(\tau)
    \end{pmatrix}, 
    \quad    
    \mathbf{\Tilde{F}}(\tau) =
    \begin{pmatrix}
        \Tilde{F}(\tau)\\
        \Tilde{G}(\tau)
    \end{pmatrix}
\end{equation}
and $\delta \big( \bm{\sigma} - \zeta \big)$ is a shorthand for a vector whose elements are $\delta( \sigma_i - \zeta)$ (\textit{i.e.} a list of $\delta$ functions that ``turn on'' the coefficients whenever the particle worldline $\zeta(\tau)$ crosses a grid point $\sigma_i$). Here, the product 
$\mathbf{\Tilde{F}}(\tau) \delta \big( \bm{\sigma} - \zeta \big) =\mathbf{\Tilde{F}}(\tau)\circ \delta \big( \bm{\sigma} - \zeta \big)$ is the Hadamard (or element-wise) product of two vectors, while  $\mathbf{L} \,\mathbf{u}$ is the inner (or dot) product of a matrix and a vector.

We now pose Eq. \eqref{eq:ch6GeneralMatrixEqn} as an integral equation over an interval $[\tau_\nu,\tau_{\nu+1}]$:
\begin{equation}
    \mathbf{u}^{\nu+1} = \mathbf{u}^\nu + \mathbf{L}  \int_{\tau_\nu}^{\tau_{\nu+1}} \mathbf{u}(\tau)~ d\tau + \int_{\tau_\nu}^{\tau_{\nu+1}} \mathbf{s}(\tau) d\tau + \int_{\tau_\nu}^{\tau_{\nu+1}} \mathbf{\Tilde{F}}(\tau) \delta( \bm{\sigma} - \zeta)d\tau
\end{equation}
and demonstrate that  discontinuous Hermite integration  may be applied to numerically solve this equation for both Solutions I and II. 
\subsubsection{Discontinuous trapezium rule}
We first apply the discontinuous trapezium rule DH2 (Eq. \eqref{eq:discTrap}) to approximate the above integrals, so
\begin{equation}
    \mathbf{u}^{\nu+1} = \mathbf{u}^\nu +\frac{\Delta \tau}{2} \mathbf{L}  (\mathbf{u}^{\nu}+\mathbf{u}^{\nu+1}) + \mathbf{K}_{\rm H2}(\Delta \tau, \Delta \tau_i) + \frac{\Delta \tau}{2}( \mathbf{s}^{\nu} + \mathbf{s}^{\nu+1}) + \llbracket \mathbf{u} \rrbracket
\end{equation}
where
\begin{equation}
    \llbracket \mathbf{u} \rrbracket = \int_{\tau_\nu}^{\tau_{\nu+1}} \mathbf{\Tilde{F}}(\tau) \delta( \bm{\sigma} - \xi)d\tau = \frac{1}{|d\zeta/d\tau|_{\tau_i}} \mathbf{\Tilde{F}}(\tau_i)~ \theta(\tau_{\nu+1} - \tau_i)~ \theta( \tau_i - \tau_\nu) 
\end{equation}
is a function that turns on when the worldline $\zeta(\tau)$ crosses the gridpoint $\sigma_i$ to alter the value of $u_i$ and $\mathbf{K}_{\rm H2}(\Delta \tau,\Delta \tau_i)$ is a vector including the jumps in $\mathbf{u}$ which appear in Eq. \eqref{eq:discTrap}; $\Delta \tau_i=\tau_i-\tau_\nu$ represents the interval from $\tau_\nu$ to the crossing time $\tau_i$ satisfying $\zeta(\tau_i)=\sigma_i$ for some $i$, if such a time exists in the interval $[\tau_\nu,\tau_{\nu+1}]$. This algebraic equation may be solved for $\mathbf{u}^{\nu+1}=\mathbf{u}(\tau_{\nu+1})$ using the methods of \cite{Markakis:2014nja} to arrive at a form which mitigates round-off error:
\begin{equation}
    \mathbf{u}^{\nu+1} = \mathbf{u}^\nu + \bigg( \mathbf{I} - \frac{\Delta \tau}{2} \mathbf{L} \bigg)^{-1}  \bigg( \Delta \tau~ \mathbf{L} \; \mathbf{u}^\nu + \mathbf{K}_{\rm H2}(\Delta \tau, \Delta \tau_i) + \frac{\Delta \tau}{2} ( \mathbf{s}^{\nu} + \mathbf{s}^{\nu+1}) + \llbracket \mathbf{u} \rrbracket \bigg)
\end{equation}
As discussed in Section \ref{sec:ch5DiscMethods}, such a scheme should exhibit second order convergence in $\Delta \tau$.
\subsubsection{Discontinuous Hermite rule}
We next apply the discontinuous Hermite rule DH4 (Eq. \eqref{eq:discHerm}) to the above integrals and find that
\begin{multline}
    \mathbf{u}^{\nu+1} = \mathbf{u}^\nu +\frac{\Delta \tau}{2} \mathbf{L}  (\mathbf{u}^{\nu}+\mathbf{u}^{\nu+1}) +\frac{\Delta \tau^2}{12} \mathbf{L}  (\Dot{\mathbf{u}}^{\nu}-\Dot{\mathbf{u}}^{\nu+1}) + \mathbf{K}_{\rm H4}(\Delta \tau, \Delta \tau_i)\\
    + \frac{\Delta \tau}{2}( \mathbf{s}^{\nu} + \mathbf{s}^{\nu+1}) + \frac{\Delta \tau^2}{12}( \Dot{\mathbf{s}}^{\nu} - \Dot{\mathbf{s}}^{\nu+1}) + \llbracket \mathbf{u} \rrbracket
\end{multline}
where now $\mathbf{K}_{\rm H4}(\Delta \tau,\Delta \tau_i)$ is a vector including the jumps in $\mathbf{u}$ which appear in Eq.~\eqref{eq:discHerm}. The $\tau-$derivatives of $\mathbf{u}$ represented by the overdot may be removed by using the original evolution equation Eq. \eqref{eq:ch6GeneralMatrixEqn}. The resulting algebraic equation may be solved for $\mathbf{u}^{\nu+1}$ using the same method as \cite{Markakis:2014nja}:
\begin{multline}
    \mathbf{u}^{\nu+1} = \mathbf{u}^\nu + \bigg( \mathbf{I} - \frac{\Delta \tau}{2} \mathbf{L}  \bigg( \mathbf{I} - \frac{\Delta \tau}{6} \mathbf{L} \bigg) \bigg)^{-1}  \bigg( \Delta \tau~ \mathbf{L}  \bigg( \mathbf{u}^\nu + \frac{\Delta \tau}{12} ( \mathbf{s}^{\nu} - \mathbf{s}^{\nu+1}) \bigg)\\
    + \mathbf{K}_{\rm H4}(\Delta \tau, \Delta \tau_i) + \frac{\Delta \tau}{2} ( \mathbf{s}^{\nu} + \mathbf{s}^{\nu+1}) + \frac{\Delta \tau^2}{12}( \Dot{\mathbf{s}}^{\nu} - \Dot{\mathbf{s}}^{\nu+1}) + \llbracket \mathbf{u} \rrbracket \bigg)
\end{multline}
In the above formula, terms polynomial in $\Delta \tau$ have been factored in Horner form to improve arithmetic precision.

We've replaced matrix multiplications with separate matrix multiplication and addition steps to reduce round-off error accumulation, as shown by prior work \cite{markakis2019timesymmetry,moboylethesis}. This enhances numerical energy and phase-space volume conservation. Modern CPU and GPU libraries accelerate these operations, allowing efficient implementation without intensive programming. Given their conservation qualities, reduced errors, and library support, these time-symmetric schemes excel in long-term evolution of distributionally sourced PDEs, like those in black hole perturbation theory and Extreme Mass Ratio Inspiral simulations. We now present numerical results.

\subsubsection{Numerical results}
With these two integration rules, we may now numerically integrate Eq.~\eqref{eq:DiscFirstOrder}
in a hyperboloidal coordinate chart and compare to
Solutions I and II  given by Eqs.~\eqref{eq:ch6FieldIHyper} and \eqref{eq:ch6FieldIIHyper} respectively. We take as initial data each solution evaluated at $\tau_0$ satisfying $t(\tau_0,\xi(\tau_0)) = 0$. 

Our schemes work for both finite-difference and pseudo-spectral methods. 
For instance, interpolation on Chebyshev-Gauss-Lobatto nodes
\begin{equation}  \label{ChebyshevNodes}
\sigma_i=\frac{{{a} + {b}}}{2} + \frac{{{b} -a}}{2}z_i, 
\quad z_i=-\cos \theta_i,
\quad
\theta_i= {\frac{{i\pi }}{N}}, \quad i=0,1,\dots,N
\end{equation}
converges uniformly 
for every absolutely continuous function \cite{Fornberg1998}. If all the function derivatives are bounded on the  interval $x \in [a,b]$,
Chebyshev interpolation has the property of exponential convergence  on the interval $x \in [a,b]$. In Section~\eqref{sec:DiscColl}, we generalized this concept to piecewise smooth functions. 
The elements $D_{ij}$ of  the $(N+1)\times (N+1)$ Chebyshev first derivative
matrix $D$ are given by
\begin{equation}
    {D}^{(1)}_{i j} = \frac{2}{b-a}    
    \begin{cases}
        \frac{c_i (-1)^{i+j}}{c_j (z_i-z_j)} & i \neq j\\
        \frac{z_i}{2(1-z_i^2)} & i = j \neq 0,N\\
        -\frac{2 N^2 + 1}{6} & i = j = 0\\
        \frac{2 N^2 + 1}{6} & i = j = N
    \end{cases}
\end{equation}
where $c_0 = c_N = 2$ and $c_1, \dotsc, c_{N-1} = 1$. The second derivative operator can be evaluated by $\mathbf{D}^{(2)} = (\mathbf{D}^{(1)} )^2$ or, equivalently \cite{canuto_2006},
\begin{equation}\label{ChebyshevD2}
    {D}^{(2)}_{i j} = \left( \frac{2}{b-a}  \right)^2
    \begin{cases}
    \frac{(-1)^{i+j}}{c_j} 
        \frac{z_i^2+z_i z_j-2}{(1-z_i^2) (z_i-z_j)^2} & i   \neq 0, \ i \neq N, \ i \neq j\\
        -\frac{(N^2-1)(1-z_i^2)+3}{3(1-z_i^2)^2}  &  i \neq 0, \ i \neq N, \ i = j\\
       \frac{2}{3} \frac{(-1)^{i}}{c_i} \frac{(2N^2+1)(1-z_i)-6}{(1-z_i)^2} & i = 0, \ j \ne 0\\
        \frac{2}{3} \frac{(-1)^{i+N}}{c_i} \frac{(2N^2+1)(1+z_i)-6}{(1+z_i)^2} & i = N, \ j \ne N\\
        -\frac{N^4 - 1}{15} & i = j = 0, \ i=j=N
    \end{cases}.
\end{equation}

\begin{figure}[t]
    \centering
    \begin{subfigure}{0.45\textwidth}
        \centering
        \includegraphics[width=\textwidth]{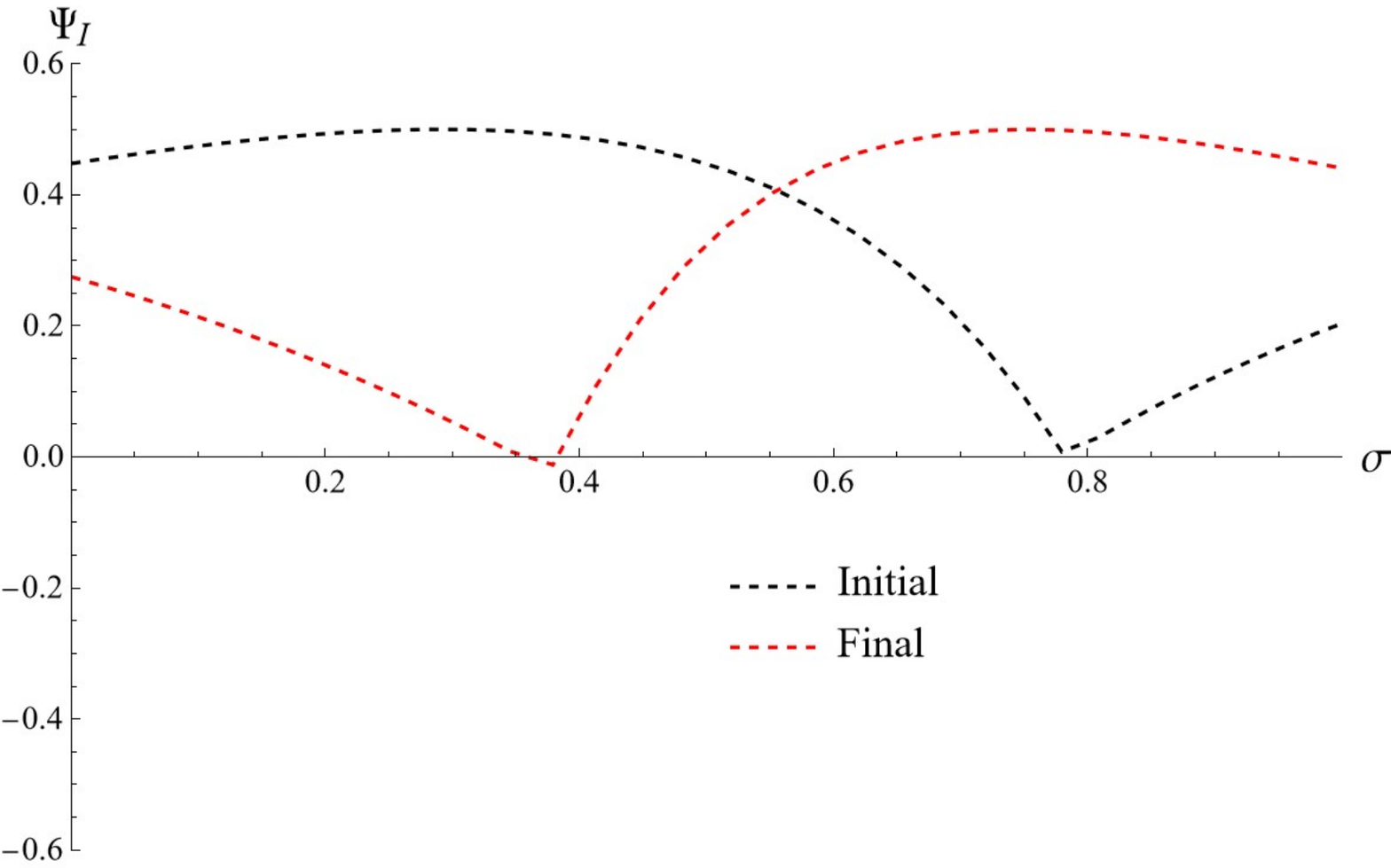}
        \caption{Solution I}
    \end{subfigure}
    \hfill
    \begin{subfigure}{0.45\textwidth}
        \centering
        \includegraphics[width=\textwidth]{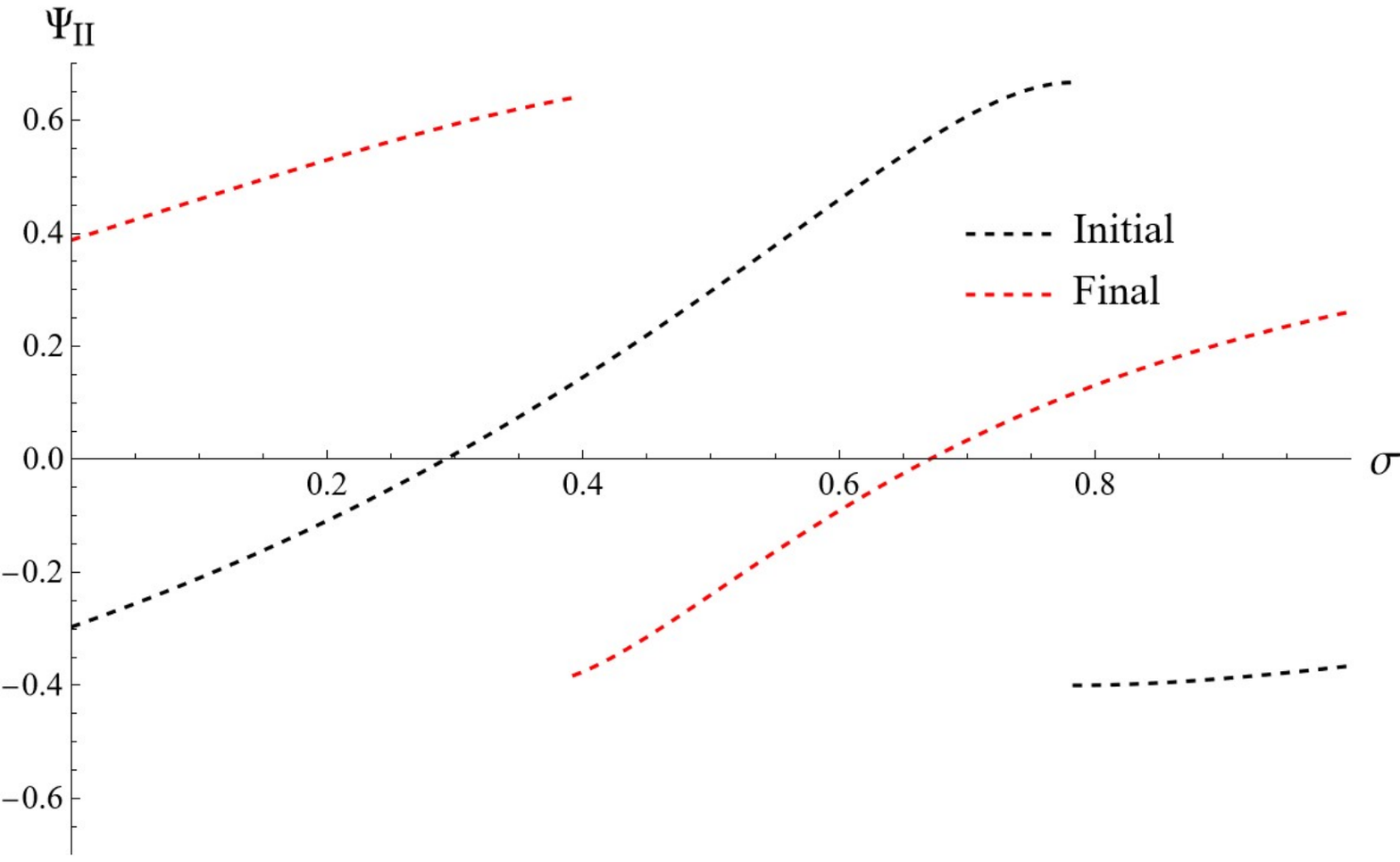}
        \caption{Solution II}
    \end{subfigure}
    \caption{Numerical solutions to Eq. \eqref{eq:ch6FullFlatEqn} using discontinuous collocation methods for spatial discretization and the discontinuous trapezoidal rule in Eq. \eqref{eq:discTrap} for time integration. The Initial function was selected to be the exact solutions provided in \cite{Field2009} at $\tau=\tau_0$; this ensures the numerical solution can be directly compared to an exact solution when studying convergence. We evolved the reduced system from $\tau=\tau_0$ to $\tau=4.3$ taking $\Delta \tau = 0.05$.}
    \label{fig:ch6FieldSol}
\end{figure}

\begin{figure}[!t]
    \centering
    \begin{subfigure}{0.45\textwidth}
        \includegraphics[width=\textwidth]{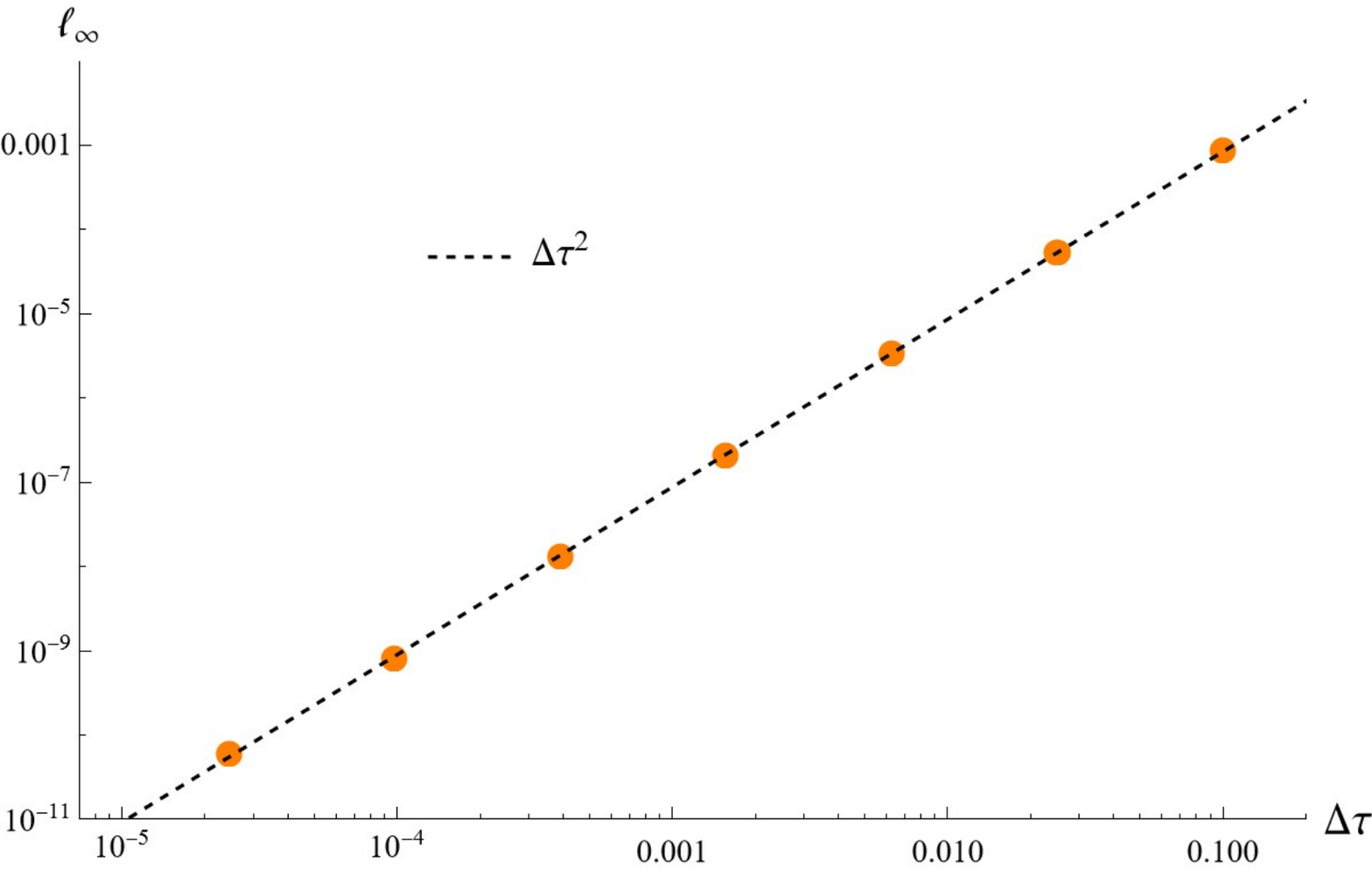}
        \caption{DH2 Temporal Convergence, Solution I}
    \end{subfigure}
    \hfill
    \begin{subfigure}{0.45\textwidth}
        \includegraphics[width=\textwidth]{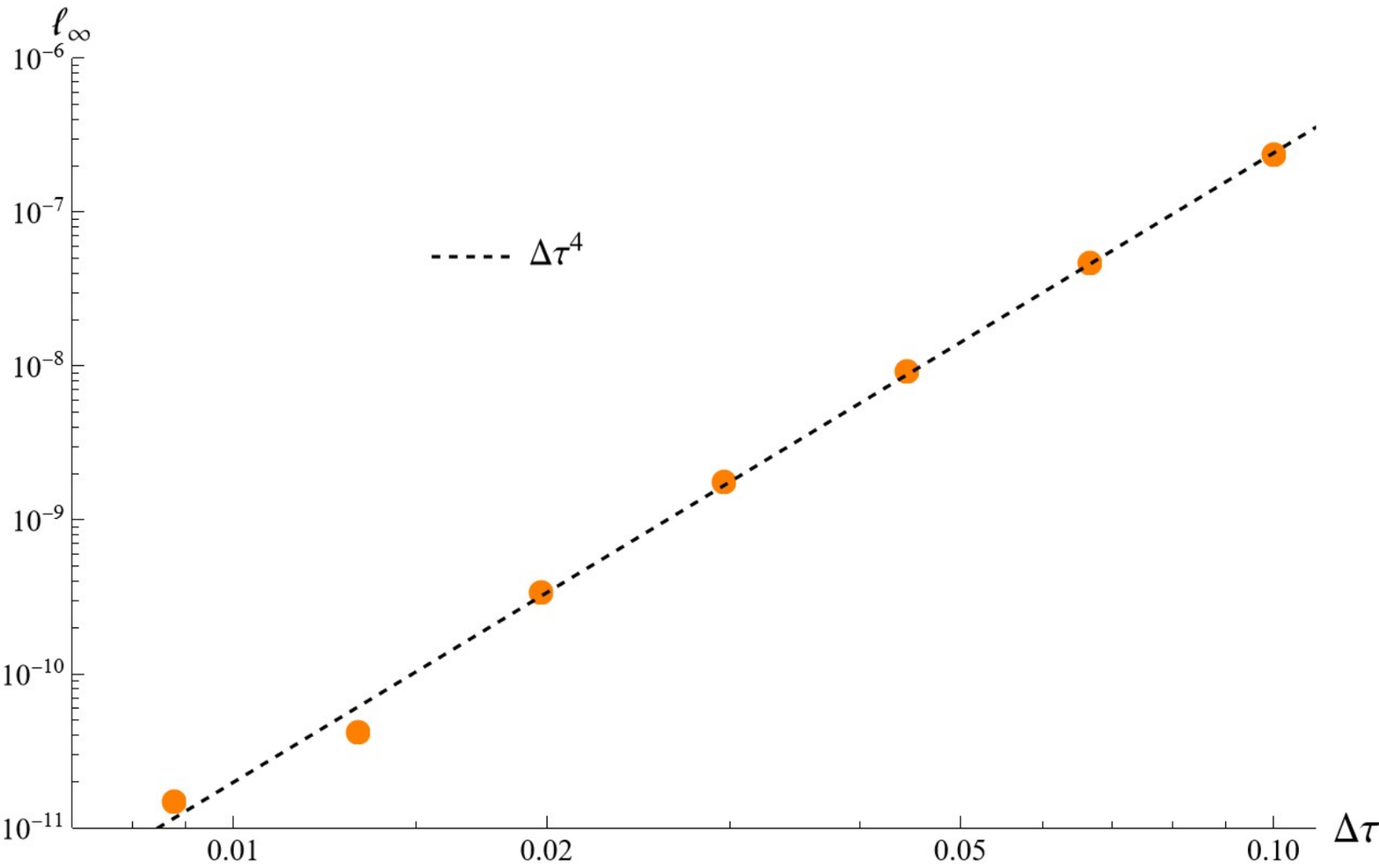}
        \caption{DH4 Temporal Convergence, Solution I}
    \end{subfigure}
    \vfill
    \begin{subfigure}{0.45\textwidth}
        \includegraphics[width=\textwidth]{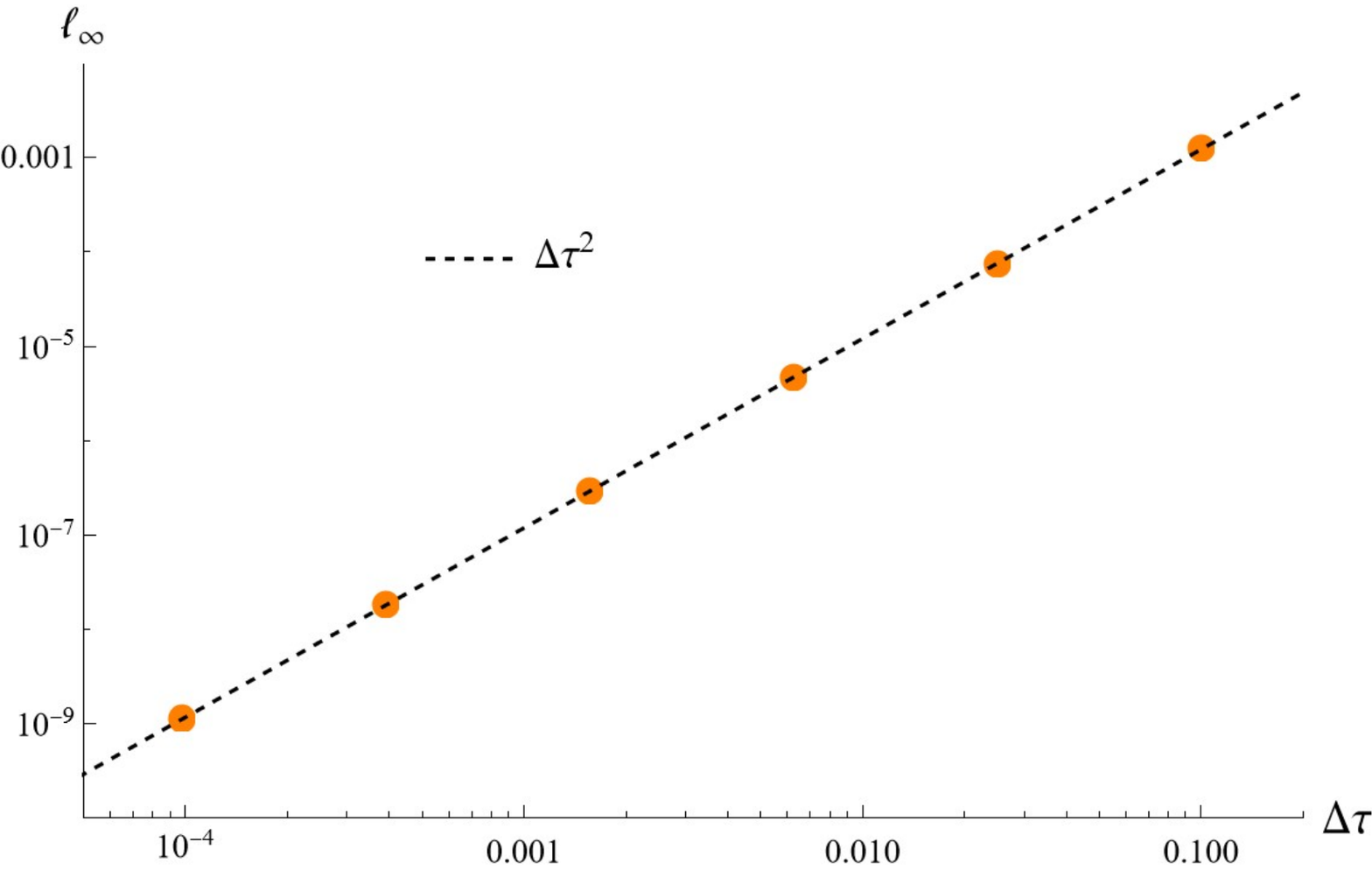}
        \caption{DH2 Temporal Convergence, Solution II}
    \end{subfigure}
    \hfill
    \begin{subfigure}{0.45\textwidth}
        \includegraphics[width=\textwidth]{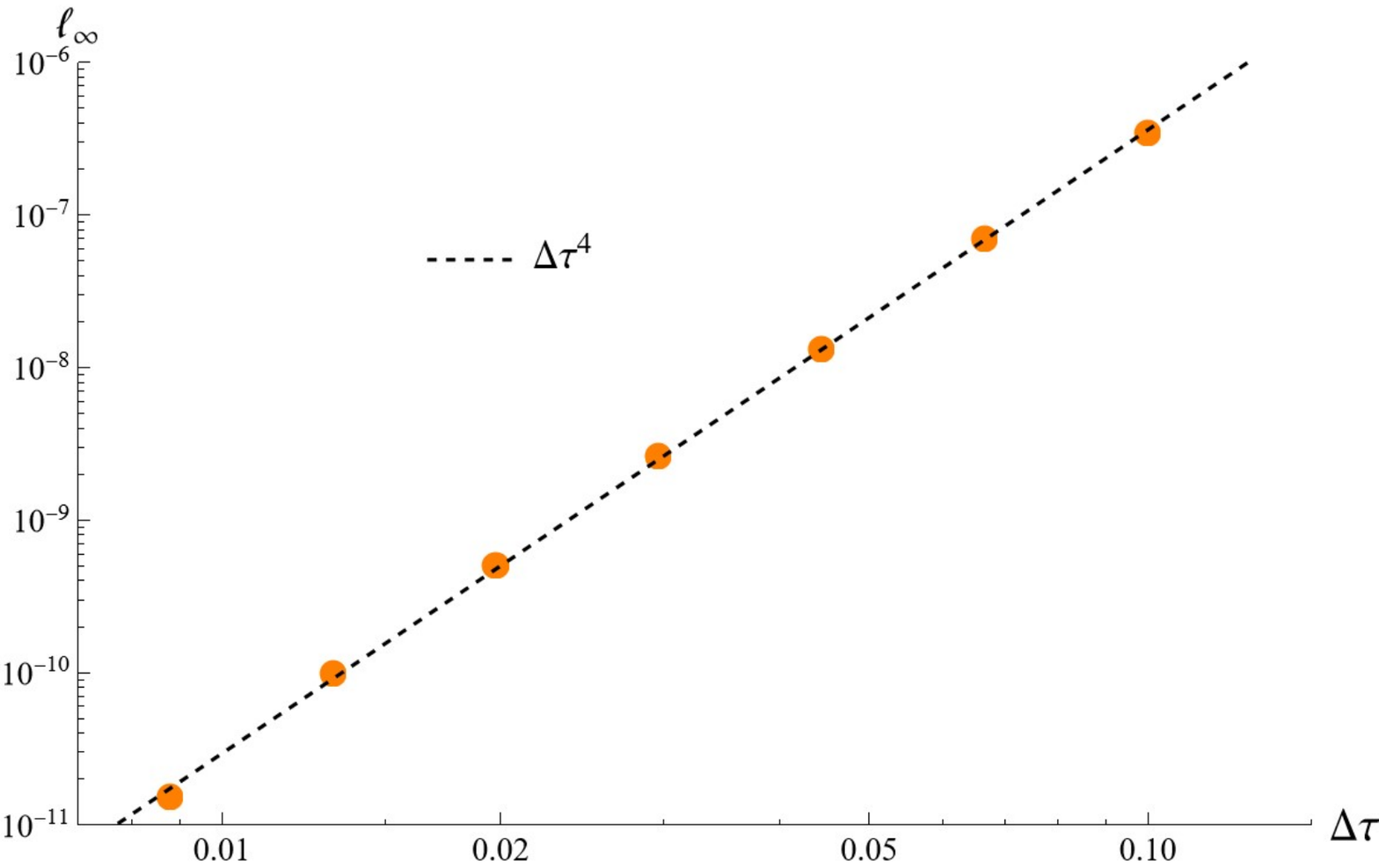}
        \caption{DH4 Temporal Convergence, Solution II}
    \end{subfigure}
    \caption{The $l_\infty$ error norm of the numerical approximation to Solutions I (Eq. \eqref{eq:ch6FieldIHyper}) and II (Eq. \eqref{eq:ch6FieldIIHyper}) using discontinuous collocation and numerical integration. In both cases, the discontinuous DH2 method error scales as $\Delta \tau^2$ and the discontinuous DH4 method error scales as $\Delta \tau^4$ , as expected.}
    \label{fig:ch6FlatHError}
\end{figure}

Here, we use a Chebyshev pseudo-spectral method (Eqs. \eqref{ChebyshevNodes} - \eqref{ChebyshevD2}) to discretize the  spatial interval $\sigma \in [0,1]$ with $N=45$ nodes and use $\simeq 10-12$ jumps for each solution. 
The results for the discontinuous trapezium rule when $\Delta \tau = 0.05$ are plotted in Figure~\ref{fig:ch6FieldSol}. 
Moreover, we observe the expected convergence for each discontinuous method as illustrated in Figure~\ref{fig:ch6FlatHError}.
 
%}

\newpage
%\\\\

%
%
%\begin{equation}\label{eq:LD4OperatorPrecise}
% \mathbf{u}(\tau_{n+1})  = \mathbf{u}(\tau_n) +{\left[ {{\mathbf{I}} - \frac{{\mathbf{L}\Delta \tau }}{2} \mkern-5mu \cdot \mkern-5mu \left( {{\mathbf{I}} - \frac{{\mathbf{L}\Delta \tau }}{5} \mkern-5mu \cdot \mkern-5mu  \left( {{\mathbf{I}} - \frac{{{\mathbf{L}}\Delta \tau }}{{{\text{12}}}}} \right)} \right)} \right]^{ - 1}} \mkern-5mu \cdot {\mathbf{L}}\Delta \tau  \cdot \left( {{\mathbf{I}} + \frac{{\mathbf{L}\Delta \tau }}{60} } \right)  \cdot \mathbf{u}(\tau_n)
%\end{equation}

\newpage
\section*{Acknowledgments}
We thank Charles F. Gammie for co-supervising the first author's PhD dissertation \cite{moboylethesis}, where the present work originally appeared. 
We thank Leor Barack for discontinuous second-order finite-differencing schemes and systematized spatial derivative jump calculations that inspired this work. We thank Pablo Brubeck for notes on using the Frobenius method to compute recursion relations between spatial derivative jumps, and for noticing that the original discontinuous collocation method converged for stationary but not for moving distributional sources; this observation motivated the creation of discontinuous time integration methods. We thank Garbielle Allen and Ed Seidel for support at NCSA. We thank John L. Friedman for valuable discussions on symplecticity and volume preserving and for suggesting a discontinuous method-of-lines graph. We thank An{\i}l Zengino\u{g}lu and Juan A. Valiente-Kroon for insights into hyperboloidal slicing and Rodrigo Panoso Macedo for suggestions and discussions on the minimal gauge. We thank Lidia Gomes Da Silva for discussions and efforts to verify our results.

\appendix

\section*{Appendix}
\section{Derivation of the recursion relation}\label{ch:recursion}
We proceed by a method of undetermined coefficients assuming the form of Eq. \eqref{eq:ch6NSFunc} with the goal of determining the $J_n(\tau)$. We first note that the properties of the unit jump functions (\textit{cf.} \cite{Markakis:2014nja}) give rise to the following identities:
\begin{equation}
    \partial_\sigma \Psi_{\mathrm{NS}} = \sum_{n=0}^\infty J_{n+1}(\tau) \Phi_{n}(\sigma;\zeta) + J_0(\tau) \delta(\sigma-\zeta)
\end{equation}
\begin{equation}
    \partial_\sigma^2 \Psi_{\mathrm{NS}} = \sum_{n=0}^\infty J_{n+2}(\tau) \Phi_n(\sigma;\zeta) + J_1(\tau) \delta(\sigma - \zeta)+ J_0(\tau) \delta'(\sigma - \zeta)
\end{equation}
\begin{equation}\label{eq:PsiTauDerivative}
    \partial_\tau \Psi_{\mathrm{NS}} = \sum_{n=0}^{\infty} \Big(\Dot{J}_n(\tau) - J_{n+1}(\tau) \Dot{\zeta} \Big)  \Phi_{n}(\sigma;\zeta) - J_0(\tau) \Dot{\zeta} \delta(\sigma-\zeta)
\end{equation}
\begin{multline}\label{eq:PsiTau2Derivative}
    \partial_\tau^2 \Psi_{\mathrm{NS}} = \sum_{n=0}^{\infty} \Big( \Ddot{J}_n(\tau) - 2 \Dot{J}_{n+1}(\tau) \Dot{\zeta} - J_{n+1}(\tau) \Ddot{\zeta} + J_{n+2}(\tau) ~\Dot{\zeta}^2 \Big) \Phi_{n}(\sigma;\zeta) + J_0(\tau) \Dot{\zeta}^2 \delta'(\sigma - \zeta)\\
    + \Big( J_1(\tau) \Dot{\zeta}^2 - J_0(\tau) \Ddot{\zeta} - 2 \Dot{J}_0(\tau) \Dot{\zeta} \Big) \delta(\sigma - \zeta)
\end{multline}
\begin{multline}
    \partial_\sigma \partial_\tau \Psi_{\mathrm{NS}} = \sum_{n=0}^{\infty} \Big(\Dot{J}_{n+1})\tau) - J_{n+2}(\tau) \Dot{\zeta} \Big)  \Phi_{n}(\sigma;\zeta) + \Big(\Dot{J}_0(\tau) - J_1(\tau) \Dot{\zeta} \Big) \delta(\sigma - \zeta)\\
    - J_0(\tau) \Dot{\zeta}~ \delta'(\sigma-\zeta)
\end{multline}
Upon substitution into the left hand side of Eq. \eqref{eq:ch6GeneralWave}, the ``product rule'' of \cite{moboylethesis} allows each term to be transformed into a pure expansion in jump functions.
\begin{equation}
    E(\sigma) \partial_\sigma \Psi_{\mathrm{NS}} = \sum_{n=0}^\infty \Bigg( \sum_{l=0}^n \binom{n}{l} E^{(l)}(\zeta)J_{n+1-l}(\tau) \Bigg) \Phi_{n}(\sigma;\zeta) + E(\zeta) J_0(\tau)\delta(\sigma-\zeta)
\end{equation}
\begin{multline}
    C(\sigma) \partial_\sigma^2 \Psi_{\mathrm{NS}} = \sum_{n=0}^\infty \Bigg( \sum_{l=0}^n \binom{n}{l} C^{(l)}(\zeta) J_{n+2-l}(\tau) \Bigg) \Phi_{n}(\tau,\sigma)\\
    + \Big(C(\zeta) J_1(\tau) - C'(\zeta) J_0(\tau) \Big) \delta(x-\zeta) + C(\zeta) J_0(\tau) \delta'(\sigma-\zeta)
\end{multline}
\begin{multline}
    B(\sigma) \partial_\tau \Psi_{\mathrm{NS}} = \sum_{n=0}^{\infty} \Bigg( \sum_{l=0}^n \binom{n}{l} B^{(l)}(\zeta) \Big(\Dot{J}_{n-l}(\tau) - \Dot{\zeta} J_{n+1-l}(\tau) \Big) \Bigg) \Phi_{n}(\sigma,\zeta)\\
    - B(\zeta) J_0(\tau) \Dot{\zeta} \delta(\sigma-\zeta)
\end{multline}
\begin{multline}
    Z(\sigma) \partial_\tau^2 \Psi_{\mathrm{NS}}\\
    = \sum_{n=0}^\infty \Bigg(\sum_{l=0}^n \binom{n}{l} Z^{(l)}(\zeta) \Big(\Ddot{J}_{n-l}(\tau) - 2 \Dot{J}_{n+1-l}(\tau) \Dot{\zeta} - J_{n+1-l}(\tau) \Ddot{\zeta} + J_{n+2-l}(\tau) \Dot{\zeta}^2 \Big) \Bigg )\Phi_n(\sigma;\zeta)\\
    + Z(\zeta) J_0(\tau) \Dot{\zeta}^2 \delta'(x - \zeta) + \Big( Z(\zeta) ( J_1(\tau) \Dot{\zeta}^2 - J_0(\tau) \Ddot{\zeta} - 2 \Dot{J}_0(\tau) \Dot{\zeta} )\\
    - Z'(\zeta) J_0(\tau) \Dot{\zeta}^2 \Big) \delta(\sigma - \zeta)
\end{multline}
\begin{multline}
    A(\sigma) \partial_\tau \partial_\sigma \Psi_{\mathrm{NS}} = \sum_{n=0}^{\infty} \Bigg( \sum_{l=0}^n \binom{n}{l} A^{(l)}(\zeta) \Big(\Dot{J}_{n+1-l}(\tau) - J_{n+2-l}(\tau) \Dot{\zeta} \Big) \Bigg)  \Phi_{n}(\sigma;\zeta)\\
    + \Big( A(\zeta) (\Dot{J}_0(\tau) - J_{1} (\tau)\Dot{\zeta} ) - A'(\zeta) J_0(\tau) \Big) \delta(x-\zeta) - J_0(\tau) A(\zeta) \Dot{\zeta} \delta'(\sigma-\zeta)
\end{multline}
These expressions are then substituted into the left hand side of Eq. \eqref{eq:ch6GeneralWave}. Equating the coefficients of $\delta'(\sigma-\zeta)$ and $\delta(\sigma-\zeta)$ results in Eqs. \eqref{eq:ch6NewJ0} and \eqref{eq:ch6NewJ1}, respectively. To obtain the recursion relation of Eq. \eqref{eq:ch6NewRecur}, observe that the coefficients of each remaining $\Phi_n$ must vanish and identify the terms proportional to $J_{n+2}(\tau)$ in each sum over $l$.

\section{Jump expressions}\label{ch:Jumps}
We present the first few jumps of Solutions I and II when transformed to the hyperboloidal coordinate chart of Eqs. \eqref{eq:compact_g} and \eqref{eq:height_h}. This expressions follow from generalized jump conditions of Eqs. \eqref{eq:ch6NewJ0}, \eqref{eq:ch6NewJ1}, and \eqref{eq:ch6NewRecur} along with the transformed source terms of Eqs. \eqref{eq:ch6NewDelta} and \eqref{eq:ch6NewDeltaP}. The first few jumps of Solution I are
\begin{subequations}
\begin{align}
    J_0 &= 0\\
    J_1 &= -\frac{1 - v(1-2\zeta^2)}{2(1-v^2)(1-\zeta)\zeta^2} \cos t_\zeta\\
    J_2 &= \frac{1}{2} \bigg( \frac{2 + \zeta}{(1+v) \zeta^3} - \frac{1}{(1-v)(1-\zeta)^2} \bigg) \cos t_\zeta - \frac{(1 - v - 2\zeta^2)(1 - v + 2 \zeta^2)}{2( 1- v^2)^2 (1-\zeta)^2 \zeta^4} \sin t_\zeta
\end{align}
\begin{multline}
    J_3 = \frac{- 3 (1-v)^3 \left(4 v^2+8 v+3\right) \zeta^4 - 3 (1-v)^3 \left(2 v^2+4 v+3\right) \zeta^2+\left(-4 v^5+10 v^3+2 v\right) \zeta^6}{2 (1-v)^3 (1+v)^3 (1-\zeta)^3 \zeta^6}\\
    \frac{+16 (1-v)^3 (1+v)^2 \zeta^3+(1-v)^3}{} \cos t_\zeta + \frac{3}{2} \bigg( \frac{(1+\zeta)(2+\zeta)}{(1+v)^2\zeta^5} + \frac{1}{(1-v)^2(1-\zeta)^3} \bigg) \sin t_\zeta
\end{multline}
\begin{multline}
    J_4 = 3 \bigg( -\frac{2+\zeta(5-v(v+2)\zeta(4+\zeta))}{(1+v)^2\zeta^5} + \frac{v(2-v)}{(1-v)^3(1-\zeta)^4} \bigg) \cos t_\zeta\\
    + \frac{1}{2} \left(\frac{-11 (2-v) v+10}{(1-v)^4 (1-\zeta)^4}+\frac{\zeta (\zeta (-(11 v (v+2)+10) \zeta (\zeta +4)-6 (6 v (v+2)+5))+4)+1}{(1+v)^4 \zeta^8}\right) \sin t_\zeta
\end{multline}
\begin{multline}
    J_5 = \frac{1}{2} \bigg( \frac{-10(2-v) v (24 (2-v) v+13)}{(1-v)^5 (1-\zeta)^5}+\frac{\zeta (\zeta (\zeta (-(v (2+v) (24 v (v+2)+13)-10) \zeta (\zeta+5)}{(1+v)^5 \zeta^{10}}\\
    \frac{+260 v(2+v)+250)+120 v (2+v)+110)-5)-1}{} \bigg) \cos t_\zeta\\
    + 5 \bigg(\frac{5 (2-v) v-4}{(1-v)^4 (\zeta -1)^5} +\frac{\zeta (\zeta((5 v (v+2)+4) \zeta (\zeta +5)+3 (8 v (v+2)+5))-7)-2}{(v+1)^4 \zeta^9} \bigg) \sin t_\zeta
\end{multline}
\begin{multline}
    J_6 = \frac{15}{2} \bigg(\frac{(2-v) v (1-8 (2-v) v)-6}{(1-v)^5 (1-\zeta)^6}    +\frac{\zeta(\zeta ((v (v+2)(8 v (v+2)+1)-6) (\zeta+6) \zeta^2}{(v+1)^5 \zeta^{11}} \\
    \frac{-14 (11 v (v+2)+10) \zeta-80 v (v+2)-64)+9)+2}{}\bigg) \cos t_\zeta\\
    + \frac{1}{2} \bigg( \frac{190 - (2-v) v (463 - 274v (2-v))}{(1-v)^6 (1-\zeta)^6}
    +\frac{\zeta(\zeta (\zeta(\zeta (-(v (v+2) (274 v (v+2)+463)}{} \\
    \frac{+190) \zeta(\zeta + 6)-15 (v (v+2) (120 v (v+2)+169)+50)) +940 v (v+2)+920)}{(v+1)^6 \zeta ^{12}}\\
    \frac{+300 v (v+2)+285)-6)-1}{}\bigg) \sin t_\zeta
\end{multline}
\end{subequations}
The first few jumps of Solution II are
\begin{subequations}
\begin{align}
    J_0 &= -\frac{\cos t_\zeta}{1- v^2}\\
    J_1 &= \frac{2(1+v^2)\zeta^2 -(1-v^2)}{2(1-v^2)^2\zeta^2(1-\zeta)} \sin t_\zeta\\
    J_2 &= \frac{1}{2} \bigg( \frac{(1+\zeta)^2}{(1+v)^3 \zeta^4} + \frac{1}{(1-v)^3(1-\zeta)^2} \bigg) \cos t_\zeta + \frac{1}{2} \bigg( \frac{2+\zeta}{(1+v)^2 \zeta^3} + \frac{1}{(1-v)^2(1-\zeta)^2} \bigg) \sin t_\zeta
\end{align}
\begin{multline}
    J_3 = \frac{3}{2} \bigg(-\frac{(1+\zeta)(2+\zeta)}{(1+v)^3 \zeta^5} + \frac{1}{(1-v)^3(1-\zeta)^3)} \bigg) \cos t_\zeta \\
    + \frac{1}{2} \bigg( \frac{1 -2v (2-v)}{(1-v)^4 (1-\zeta)^3} + \frac{\zeta (3-(2 v (v+2)+1) \zeta (\zeta +3)) +1}{(v+1)^4 \zeta^6} \bigg) \sin t_\zeta
\end{multline}
\begin{multline}
    J_4 = \frac{1}{2} \bigg(\frac{10-11 v (2-v)}{(1-v)^5 (1-\zeta)^4} + \frac{\zeta (\zeta ((11 v (v+2)+10) \zeta(\zeta +4)+6 (6 v (v+2)+5))-4)-1}{(v+1)^5 \zeta^8} \bigg) \cos t_\zeta\\
    + 3 \bigg(-\frac{v(2-v)}{(1-v)^4 (1-\zeta)^4} + \frac{\zeta (v (v+2) \zeta (\zeta +4)-5)-2}{(v+1)^4 \zeta^7} \bigg) \sin t_\zeta
\end{multline}
\begin{multline}
    J_5 = 5 \bigg(-\frac{5 (2-v) v+4}{(1-v)^5 (1-\zeta)^5} + \frac{\zeta (\zeta (-(5 v (v+2)+4) \zeta (\zeta +5)-3 (8 v (v+2)+5))+7)+2}{(v+1)^5 \zeta^9}\bigg) \cos t_\zeta \\
    + \frac{1}{2} \bigg(\frac{-10+ v(2-v) (24v (2-v) v-13)}{(1-v)^6 (1-\zeta)^5} + \frac{\zeta (\zeta (\zeta (-(v (v+2) (24 v (v+2)+13)-10) \zeta (\zeta +5)}{(v+1)^6 \zeta^{10}}\\
    \frac{+260 v (v+2)+250)+120 v (v+2)+110)-5)-1}{}\bigg) \sin t_\zeta
\end{multline}
\begin{multline}
    J_6 = \frac{1}{2} \bigg( \frac{\zeta (\zeta (\zeta ((v (v+2) (274 v (v+2)+463)+190) (\zeta +6) \zeta^2 +15 (v (v+2) (120 v (v+2)+169)}{(v+1)^7 \zeta^{12}}\\
    \frac{+50) \zeta -940 v (v+2)-920)-15 (20 v (v+2)+19))+6)+1}{}\\
    -\frac{v (2-v)(-274 v (2-v) +463)+190}{(1-v)^7 (1- \zeta)^6} \bigg) \cos t_\zeta + \frac{15}{2} \bigg(\frac{(v-2) v (8 (v-2) v+1)-6}{(1-v)^6 (1-\zeta)^6} \\
    +\frac{\zeta (\zeta ((v (v+2) (8 v (v+2)+1)-6) (\zeta+6) \zeta^2-14 (11 v (v+2)+10) \zeta}{(v+1)^6 \zeta^{11}}\\
    \frac{-80 v(v+2)-64)+9)+2}{}\bigg) \sin t_\zeta
\end{multline}
\end{subequations}

\section{Time Jumps}
As shown in Eqs.~\eqref{eq:space_jump}-\eqref{eq:time_jump}, these same jumps may be used to determine the $\tau$-direction jumps needed for the discontinuous Hermite integration schemes up to a sign. We already showed that
\begin{equation}
    K_0(\zeta) = - J_0(T).
\end{equation}
where $T = \zeta^{-1}(\sigma)$ as before. We go further by noting that Eq.~\eqref{eq:PsiTauDerivative} implies
\begin{equation}
    K_1(\zeta) = -\big(\Dot{J}_0(T) - \Dot{\zeta} J_1(T) \big)
\end{equation}
and Eq.~\eqref{eq:PsiTau2Derivative} implies
\begin{equation}
    K_2(\zeta) = -\big( \Ddot{J}_0(T) - 2 \Dot{J}_{1}(T) \Dot{\zeta} - J_{1}(T) \Ddot{\zeta} + J_{2}(T) ~\Dot{\zeta}^2 \big)
\end{equation}
To determine $K_3$, differentiate Eq.~\eqref{eq:PsiTau2Derivative} and negate the coefficient of $\Phi_0$ to find
\begin{equation}
    K_3(\zeta) = -\big( \dddot J_0(T) - 3 \Ddot{J}_1(T) \Dot{\zeta} - 3 \Dot{J}_1(T) \Ddot{\zeta} - J_1(T) \dddot \zeta + 3 \Dot{J}_2(T) \Dot{\zeta}^2 + 3 J_2(T) \Dot{\zeta} \Ddot{\zeta} - J_3(T) \Dot{\zeta}^3 \big)
\end{equation}

\newpage

\bibliography{references}

%% not want to use model5-names.bst.


        
        %% \bibitem must have one of the following forms:
        %%   \bibitem[Jones et al.(1990)]{key}...
        %%   \bibitem[Jones et al.(1990)Jones, Baker, and Williams]{key}...
        %%   \bibitem[Jones et al., 1990]{key}...
        %%   \bibitem[\protect\citeauthoryear{Jones, Baker, and Williams}{Jones
                %%       et al.}{1990}]{key}...
        %%   \bibitem[\protect\citeauthoryear{Jones et al.}{1990}]{key}...
        %%   \bibitem[\protect\astroncite{Jones et al.}{1990}]{key}...
        %%   \bibitem[\protect\citename{Jones et al., }1990]{key}...
        %%   \harvarditem[Jones et al.]{Jones, Baker, and Williams}{1990}{key}...
        %%
        
        
        % \end{thebibliography}

\end{document}